\documentclass[11pt]{amsart}

\usepackage{amsmath, amssymb, amsthm, amsrefs}
\usepackage{latexsym}
\usepackage{indentfirst}
\usepackage{graphicx}
\usepackage{placeins}
\usepackage{booktabs}
\usepackage{algorithm}
\usepackage{algorithmic}
\usepackage{multirow}

\theoremstyle{plain}

\theoremstyle{definition}

\theoremstyle{remark}

\newcommand{\bbm}{\begin{bmatrix}}
\newcommand{\ebm}{\end{bmatrix}}

\newcommand{\R}{\mathrm{R}}

\newcommand{\p}{\partial}

\begin{document}

\title[Sparsifying Preconditioner for Lippmann-Schwinger]{Sparsifying
  Preconditioner for the Lippmann-Schwinger Equation}

\author{Lexing Ying} 

\address{
  Department of Mathematics and Institute for Computational and Mathematical Engineering,
  Stanford University,
  Stanford, CA 94305
}

\email{lexing@math.stanford.edu}

\thanks{This work was partially supported by the National Science
  Foundation under award DMS-0846501 and the U.S. Department of
  Energy’s Advanced Scientific Computing Research program under award
  DE-FC02-13ER26134/DE-SC0009409. The author thanks Lenya Ryzhik for
  providing computing resources and Anil Damle for comments and
  suggestions.}



\begin{abstract}
  The Lippmann-Schwinger equation is an integral equation formulation
  for acoustic and electromagnetic scattering from an inhomogeneous
  media and quantum scattering from a localized potential. We present
  the sparsifying preconditioner for accelerating the iterative
  solution of the Lippmann-Schwinger equation. This new preconditioner
  transforms the discretized Lippmann-Schwinger equation into sparse
  form and leverages the efficient sparse linear algebra algorithms
  for computing an approximate inverse. This preconditioner is
  efficient and easy to implement. When combined with standard
  iterative methods, it results in almost frequency-independent
  iteration counts. We provide 2D and 3D numerical results to
  demonstrate the effectiveness of this new preconditioner.
\end{abstract}


\maketitle

\section{Introduction}

This paper is concerned with the efficient solution of the
Lippmann-Schwinger equation, which describes time-harmonic scattering
from inhomogeneous media in acoustics and electromagnetics as well as
time-harmonic scattering from localized potentials in quantum
mechanics. The simplest form of this equation comes from inhomogeneous
acoustic scattering. Let $\omega$ be the frequency of the
time-harmonic wave and denote the index of refraction by $1-m(x)$. The
inhomogeneity $m(x)$ is a function supported in a compact domain
$\Omega\subset \R^d$ of size $O(1)$, so the index of refraction is $1$
outside $\Omega$. Given an incoming wave $u_I(x)$ that satisfies the
free space Helmholtz equation
\[
(\Delta + \omega^2) u_I(x) = 0,\quad x\in\R^d
\]
the problem of scattering from inhomogeneous media is to find the
scattered field $u(x)$ such that the total field $u(x)+u_I(x)$
satisfies
\begin{equation}
  ( \Delta + \omega^2(1-m(x)) ) (u(x)+u_I(x)) = 0, \quad x\in\R^d
  \label{eq:helm}
\end{equation}
and $u(x)$ obeys the Sommerfeld radiation condition:
\begin{equation}
  \lim_{r\rightarrow\infty} r^{(d-1)/2} \left(\frac{\p u}{\p r} - i \omega u\right) = 0.
  \label{eq:somm}
\end{equation}
For the free space Helmholtz operator $-(\Delta +
\omega^2)$, the Green's function is given
by
\[
G(x) = 
\begin{cases}
  \frac{1}{2 i \omega} \exp(i \omega |x|), &d=1,\\
  \frac{i}{4} H^1_0(\omega |x|), & d=2,\\
  \frac{1}{4\pi |x|} \exp(i\omega |x|), & d=3.
\end{cases}
\]
Rewriting \eqref{eq:helm} as
\[
(-\Delta-\omega^2 + \omega^2 m(x)) u(x) = - \omega^2 m(x) u_I(x),\quad x\in\R^d
\]
and convolving it with $G(x)$ gives
\begin{equation}
  u + G \ast (\omega^2 m u) = G \ast (-\omega^2 m u_I),
  \label{eq:ls}
\end{equation}
which is the Lippmann-Schwinger equation written in terms of the
scattered field $u(x)$. One can also write the Lippmann-Schwinger in
terms of the total field $u(x)+u_I(x)$, but we shall stick to
\eqref{eq:ls} since its unknown $u(x)$ satisfies the Sommerfeld
radiation condition. For the inhomogeneous electromagnetic scattering
and quantum scattering, the same derivation results in integral
equations similar to \eqref{eq:ls}.

For high frequency wave fields (i.e., $\omega$ large), computing
numerical solution of scattering in inhomogeneous media remains a
challenging computational problem especially in 3D. For practical and
numerical purposes, working with \eqref{eq:ls} rather than
\eqref{eq:helm} offers several advantages.
\begin{itemize}
\item First, \eqref{eq:ls} is written in terms of $u(x)$ only for
  $x\in \Omega$ as $m(x)=0$ for $x\not\in\Omega$. Therefore, the
  unknown function in \eqref{eq:ls} is defined on the compact set
  $\Omega$ versus the whole space $\R^d$ as in \eqref{eq:helm}.
\item Second, once $u(y)$ for $y\in\Omega$ is computed, the whole
  scattered field in $\R^d$ defined via
  \[
  u(x) = - \int_\Omega G(x-y) \omega^2 m(y) (u(y)+u_I(y)) dy, \quad
  \forall x\in\R^d
  \]
  satisfies the Sommerfeld radiation condition \eqref{eq:somm}
  automatically. However, for \eqref{eq:helm} a special treatment such
  as perfectly matched layer (PML) \cite{berenger-1994,chew-1994} or
  absorbing boundary condition (ABC) \cite{engquist-77} is required to
  approximate \eqref{eq:somm} at a finite distance.
\item Third, the convolution structure of \eqref{eq:ls} allows for
  rapid application of the integral operator via the fast Fourier
  transform (FFT) \cite{cooley-1965}. 
\item Finally, most treatments of \eqref{eq:helm} for high-frequency
  problems suffer from the pollution effect \cite{babuska-1997}, due
  to the error from the local finite difference or finite element
  stencils. This problem does not show up in \eqref{eq:ls} because of
  the explicit introduction of the Green's function.
\end{itemize}

However, working with \eqref{eq:ls} also brings a couple of numerical
issues. Discretizing the integral equation \eqref{eq:ls} results a
dense linear system, which renders typical direct solvers impractical
except for 1D problems. For this reason, most numerical solutions of
the Lippmann-Schwinger equation use iterative solvers. Since the
Lippmann-Schwinger equation is numerically ill-conditioned, standard
iterative methods without preconditioning requires a huge number of
iterations for high frequency problems. This is especially true when
the inhomogeneity $m(x)$ has sharp transitions (see \cite{duan-2009}
for example). Therefore, there is a clear need for developing good
preconditioners for the Lippmann-Schwinger equation.

There has been a substantial amount of work on computing numerical
solutions of the Lippmann-Schwinger in the recent literature, for
example
\cite{andersson-2005,bruno-2004,chen-2002,lanzara-2004,sifuentes-2010,vainikko-2000}.
Several of these methods are concerned with the design of efficient
preconditioners. In \cite{bruno-2004}, Bruno and Hyde constructed a
preconditioner by replacing the inhomogeneity $m(x)$ with a piecewise
constant and radially symmetric approximation and then inverting the
associated approximate integral operator with semi-analytic
methods. Such a preconditioner is effective when there exists a good
piecewise constant and radially symmetric approximation, and naturally
deteriorates when there is no such approximation.

In \cite{chen-2002}, Chen proposed a method for solving the 2D
Lippmann-Schwinger equation using a technique that is now often
referred as recursive interpolative decomposition or recursive
skeletonization. This method generates in $O(N^{3/2})$ steps an
approximate inverse which can be used either as a direct solver if the
accuracy is sufficient or as a preconditioner. More advanced and
efficient methods for general integral equations have been proposed in
\cite{corona-2014} for 2D problems and in \cite{ho-2013} for both 2D
and 3D problems. These recent methods achieve quasi-linear complexity
for integral equations with non-oscillatory kernels, but fall back to
the same $O(N^{3/2})$ cost for 2D high frequency Lippmann-Schwinger
equations. For all these methods, the prefactor of the complexity
depends on the required accuracy and can be very large when they are
used as direct solvers.

There have also been attempts to apply two-grid and multigrid
methods to precondition the Lippmann-Schwinger equation (see
\cite{sifuentes-2010} for example). However, due to the highly
oscillatory nature of the solution field, such methods tend to improve
the iteration number at most by a constant factor.

In this paper, we introduce a new preconditioner called the
sparsifying preconditioner for the Lippmann-Schwinger equation. The
main idea is to transform the discretized Lippmann-Schwinger equation
approximately into a sparse linear system. This is possible since the
scattered field indeed satisfies the Helmholtz equation in $\Omega$ and
approximately the absorbing-type boundary conditions on
$\p\Omega$. Once the sparse linear system is ready, we invert it by
leveraging the efficiency of sparse linear algebra algorithms.

The rest of this paper is organized as follows. In Section 2 we
introduce the discretization scheme and describe the main idea.
Sections 3 and 4 explain the details of the sparsifying preconditioner
for rectangular domains and general domains, respectively. Section 5
extends our approach to the Laplace equation with
potential perturbation. Finally, discussions and future work are given in
Section 6.

\section{Discretization and main idea}

The sparsifying preconditioner to be presented is mostly independent
of the discretization scheme and the dimension. We assume that a
Nystr\"om method is used in order to keep the presentation simple.
Our discussion focuses on the 2D and 3D cases, since the 1D problem is
not computationally challenging.

At any point $x\in \Omega$, the local wave velocity is given by $c(x) =
(1-m(x))^{-1/2}$. At frequency $\omega$, the local wavelength at $x$
is equal to $\lambda(x) = (2\pi/\omega) (1-m(x))^{-1/2}$.

\subsection{Discretization}

Let us first consider the 2D case. The domain $\Omega$ is discretized
with a uniform grid of step size $h$, where $h$ is chosen such that
there are at least a few grid points per wavelength across the
domain. Such a uniform grid is convenient when the variation in $c(x)$
is not too large. For problems with large variations in $c(x)$, an
adaptively-refined grid should be used instead. The grid points are
indexed by integer pairs $i=(i_1,i_2)$ with location given by $x_i =
ih =(i_1h,i_2h)$.

We assume that $m(x)$ is zero within two layers of grid points from
the boundary $\p\Omega$, which can be easily satisfied by slightly
enlarging $\Omega$ if necessary. The discretization uses $m_i =
m(x_i)$ and let $u_i$ be the numerical approximation to
$u(x_i)$. Writing the integral on the left hand side of \eqref{eq:ls}
explicitly at the point $x_i$ gives
\begin{equation}
  (G\ast (\omega^2 m u))(x_i) = \int_\Omega G(x_i-y) (\omega^2 m(y) u(y)) dy
  \label{eq:Gmuint}
\end{equation}
and this integral is approximated with numerical quadrature using the
grid points in $\Omega$.


When $m(x)$ is discontinuous, discretizing \eqref{eq:Gmuint}
efficiently and accurately is a non-trivial task and we refer to
\cite{beylkin-2008,beylkin-2009} for recent developments. For
simplicity, we assume here that $m(x)$ is smooth (but still with
possible sharp transitions) and that the grid size $h$ is chosen to
sufficiently resolve the transitions in $m(x)$. For such $m(x)$, the
discretization of \eqref{eq:Gmuint} is more straightforward and the
only difficulty is the singularity of $G(x)$ at the origin. A simple
quadrature rule is
\[
\int_\Omega G(x_i-y) (\omega^2 m(y) u(y)) dy \approx
\sum_{j: x_j\in\Omega} k_{i-j} \omega^2 m_j u_j,
\]
where 
\[
k_t = G(ht) h^2, \quad t\not=(0,0)
\]
and $k_{(0,0)}$ is given by a quadrature correction near the
origin. This is a one-point correction of the trapezoidal rule and
gives an accuracy of O($h^4 \log(1/h))$ (see \cite{duan-2009} for
example). There are also higher order versions that correct the grid
points in the neighborhood of the origin and they do not affect the
following discussion of the sparsifying preconditioner significantly.

After quadrature approximation, the discrete form of the
Lippmann-Schwinger equation \eqref{eq:ls} takes the form, for $i$ with
$x_i\in \Omega$,
\begin{equation}
  u_i + \sum_{j:x_j\in\Omega} k_{i-j} \omega^2 m_j u_j = g_i,
  \label{eq:lsd}
\end{equation}
where $g_i$ is the discretized value of the right hand side of
\eqref{eq:ls}. In a slight abuse of notation, we define $u$ to be the
vector with entries given by $u_i$ and similarly for $g$ and
$m$. Then, \eqref{eq:lsd} can be written as
\begin{equation}
  (I + K\omega^2 m) u = g,
  \label{eq:lso}
\end{equation}
where $K$ is the matrix with entries defined via $K(i,j) = k_{i-j}$ and
$m$ is understood as operator of multiplying with the vector $m$
entry-wise.

The situation for the 3D case is almost the same except that $i$ is
now an integer triple $(i_1,i_2,i_3)$ and $k_t = G(th) h^3$ for
$t\not=(0,0,0)$.

\subsection{Main idea} 

The main idea of our approach is to find a sparse and local operator
$Q$ such that in 
\[
Q (I+ K \omega^2 m) u = Qg
\]
the operator $Q (I+ K \omega^2 m)$ is also numerically sparse and
localized. This is possible since the scattered field $u$ satisfies
\begin{itemize}
\item the Helmholtz equation \eqref{eq:helm} in the interior of $\Omega$, and 
\item approximately the absorbing-type boundary conditions on the
  boundary of $\Omega$.
\end{itemize}

For a grid point indexed by $i$, we define its neighborhood $\mu(i)$
to be
\[
\mu(i) = \{j: \|i-j\|_\infty \le 1, x_j \in \Omega \}.
\]
A grid point $i$ is called an {\em interior} point if all the grid
points in $\mu(i)$ are also in $\Omega$, and otherwise a {\em
  boundary} point. The set $\mu(i)$ contains 9 points in 2D and 27
points in 3D for an interior point, but has fewer points for a
boundary point. Let $N_I$ and $N_B$ denote the numbers of interior
and boundary grid points, respectively. The total number of grid
points $N$ is equal to $N_I + N_B$.

For the interior points, we shall design a matrix $A$ of size $N_I
\times N$, with the rows indexed by the interior grid points. For a
fixed interior index $i$, the row $A(i,:)$ satisfies two conditions:
\begin{itemize}
\item $A(i,:)$ has support in $\mu(i)$, and
\item $A(i,:) K(:,j) \equiv A(i,\mu(i)) K(\mu(i),j) \approx 0$ for any
  $j\not\in \mu(i)$, or, equivalently $A(i,:) K(:,\mu(i)^c) \approx
  0$.
\end{itemize}
Once $A$ has been calculated, we can define a sparse matrix $C$ of the
same size such that
\begin{itemize}
\item $C(i,:)$ has support in $\mu(i)$, and
\item $C(i,\mu(i)) := A(i,\mu(i)) K(\mu(i),\mu(i))$.
\end{itemize}
The above conditions imply that
\begin{itemize}
  \item $C$ and $A$ are sparse operators that represent a local
    stencil, and
  \item $A K \approx C$.
\end{itemize}
Therefore,
\begin{equation}
  A(I + K \omega^2 m) \approx (A + C \omega^2 m).
  \label{eq:A}
\end{equation}

For the boundary points, we shall design a matrix $B$ of size $N_B
\times N$, with the rows indexed by the boundary grid points. For a
fixed boundary index $i$, the row $B(i,:)$ satisfies two conditions:
\begin{itemize}
\item $B(i,:)$ has support in $\mu(i)$, and
\item $B(i,:) K(:,j) \equiv B(i,\mu(i)) K(\mu(i),j) \approx 0$ for any
  $j\not\in \mu(i)$.
\end{itemize}
Since $m_j=0$ for any grid point $j$ within two layers from the
boundary and $i$ is a boundary point, $m_j=0$ for $j\in\mu(i)$.
Therefore, the second condition implies that
\[
B(i,:) (K \omega^2 m)(:,:) \approx 0.
\]
Hence,
\begin{equation}
  B(I+ K\omega^2 m) \approx B.
  \label{eq:B}
\end{equation}

Applying matrices $A$ and $B$ to both sides of \eqref{eq:lso} yields
\begin{equation}
  \begin{bmatrix}
    A\\B
  \end{bmatrix} 
  (I+ K \omega^2 m) u = 
  \begin{bmatrix}
    A\\B
  \end{bmatrix} 
  g.
  \label{eq:AB}
\end{equation}
Combining \eqref{eq:A} and \eqref{eq:B} with \eqref{eq:AB} gives
\begin{equation}
  \begin{bmatrix}
    A + C\omega^2 m \\ B
  \end{bmatrix} 
  u \approx
  \begin{bmatrix}
    A \\ B
  \end{bmatrix}
  g.
  \label{eq:PCapp}
\end{equation}
This suggests defining the preconditioner as the mapping from $g$ to
$u$ given by
\begin{equation}
  u \Leftarrow
  \begin{bmatrix}
    A + C\omega^2 m \\ B
  \end{bmatrix}^{-1}
  \begin{bmatrix}
    A \\ B
  \end{bmatrix}
  g.
  \label{eq:PC}
\end{equation}
Since $A$ and $B$ are sparse, applying these operators on the right
hand side is fast. In addition, since $A$, $B$, and $C$ are all
designed to be local operators, the linear system solve in
\eqref{eq:PC} can be done efficiently with, for example, the nested
dissection method and the multifrontal method
\cite{duff-1983,george-1973}.


In the following sections, we shall apply this procedure to several
different setting. In each case, we explicitly give the construction
for $A$ and $B$, analyze the complexity of constructing and applying
the preconditioner, and provide numerical results.

\section{Rectangular domains}
\label{sec:R}

This section considers the case of rectangular $\Omega$. 

\subsection{Algorithm}
Let us consider the 2D case first. For simplicity, let $\Omega =
(0,1)^2$ and define $n=1/h$ to be the number of points in each
dimension. Clearly, $N=n^2$, $N_B = O(n)$, and $N_I = O(n^2)$.  In
order for the operator $B$ to serve as a sufficiently accurate
approximation to the Sommerfeld radiation condition, we assume that
$m(x)$ vanishes in a buffer region of a constant, but small, width
near the boundary and the constant $b$ denotes the ratio between the
width and the step size.

For each interior point $i$, recall that we require
\[
A(i,\mu(i)) K(\mu(i),j) \approx 0, \quad \forall j\not\in \mu(i).
\]
Since $K$ is translational invariant, it is convenient to require $A$
be so as well. By translating the point $i$ to the origin, it is
sufficient to consider the problem of finding a vector $\alpha$ such
that
\[
\alpha\cdot K(\mu(0),j) \approx 0, \quad \forall j \in I_n := \{j: -n<j_1,j_2<n, 
j\not\in \mu(0)\}.
\]
To solve for $\alpha$, we formulate the following optimization problem
\begin{equation}
  \min_{\alpha: \|\alpha\|=1} \|\alpha\cdot K(\mu(0), I_n)\|^2.
  \label{eq:Fmin}
\end{equation}
Through the singular value decomposition $ K(\mu(0),I_n) = U S V^*$,
the solution is given by
\[
\alpha = U(:,|\mu(0)|)^*
\]
and we set $A(i,\mu(i)) = \alpha$ for all interior $i$.

For the boundary, there are two separate cases: edge points and corner
points. For a boundary point $i$ on an edge (suppose the right one),
we require
\[
B(i,\mu(i)) K(\mu(i),j) \approx 0, \quad \forall j\not\in \mu(i).
\]
Translating the point $i$ to the origin gives the sufficient problem
of finding $\alpha$ such that
\[
\alpha\cdot K(\mu_E(0),j) \approx 0, \quad \forall j \in E_n := \{j: -n<j_1<-b, -n<j_2<n\},
\]
where $\mu_E(0)$ is the translated copy of $\mu(i)$. To solve for
$\alpha$, we consider
\begin{equation}
  \min_{\alpha: \|\alpha\|=1} \|\alpha\cdot K(\mu_E(0), E_n)\|^2.
  \label{eq:Emin}
\end{equation}
Once again, the singular value decomposition $ K(\mu_E(0),E_n) = U S
V^*$ gives the solution
\[
\alpha = U(:,|\mu_E(0)|)^*
\]
and we set $B(i,\mu(i)) = \alpha$ for all boundary $i$ on the edges.

For a boundary point $i$ at a corner (suppose the top-right one), we require
\[
B(i,\mu(i)) K(\mu(i),j) \approx 0, \quad \forall j\not\in \mu(i).
\]
Translating the point $i$ to the origin gives the problem of finding
$\alpha$ such that
\[
\alpha\cdot K(\mu_C(0),j) \approx 0, \quad \forall j \in C_n := \{j: -n<j_1,j_2<-b\},
\]
where $\mu_C(0)$ is the translated copy of $\mu(i)$. To solve for
$\alpha$, we consider
\begin{equation}
  \min_{\alpha: \|\alpha\|=1} \|\alpha\cdot K(\mu_C(0), C_n)\|^2.
  \label{eq:Cmin}
\end{equation}
The singular value decomposition $ K(\mu_C(0),C_n) = U S V^*$ gives the solution
\[
\alpha = U(:,|\mu_C(0)|)^*
\]
and we set $B(i,\mu(i)) = \alpha$ for all corner points $i$.

Once $A$ and $B$ are known, we compute $C$ and form
\begin{equation}
  \begin{bmatrix}
    A + C\omega^2 m \\ B
  \end{bmatrix},
  \label{eq:ABop}
\end{equation}
which is built from local compact stencils on a uniform rectangular
grid. As we mentioned earlier, the nested dissection algorithm is used
to compute an factorization of this operator.

The extension to 3D is quite straightforward, except there are three
types for boundary points: face points, edge points, corner
points. For the translated optimization problem, the point sets are
defined via
\begin{itemize}
\item interior point
  \[
  I_n := \{j: -n<j_1,j_2,j_3<n, j\not\in \mu(0)\},
  \]
\item surface point (suppose on $x_1=1$)
  \[
  F_n := \{j: -n<j_1<-b, -n<j_2,j_3<n\},
  \]
\item edge point (suppose on $x_1=1, x_2=1$)
  \[
  E_n := \{j: -n<j_1,j_2<-b, -n<j_3<n\},
  \]
\item corner point (suppose at $x_1=2,x_2=1,x_3=1$)
  \[
  C_n := \{j: -n<j_1,j_2,j_3<-b\}.
  \]
\end{itemize}

\subsection{Complexity}
Here we analyze the cost of constructing and applying the sparsifying
preconditioner.

In 2D, the setup algorithm consists of two parts: computing the
singular value decompositions and factorizing \eqref{eq:ABop} with the
nested dissection algorithm. The former has a $O(n^2)=O(N)$ cost while
the latter takes $O(n^3) = O(N^{3/2})$ steps. Therefore, the total
setup cost is $O(N^{3/2})$. Applying the preconditioner is essentially
a solve with the nested dissection algorithm, which has $O(n^2 \log n)
= O(N\log N)$ complexity.

In 3D, factorizing \eqref{eq:ABop} with the nested dissection
algorithm has an $O(n^6) = O(N^2)$ cost and this dominates the setup
cost. The cost of applying the preconditioner is equal to $O(n^4) =
O(N^{4/3})$, which is the cost of a nested dissection solve.

\subsection{Numerical results}

This preconditioner and the necessary nested dissection algorithm are
implemented in Matlab. The numerical results below are obtained on a
desktop computer with CPU speed at 2.0Hz. For the iterative solver
the GMRES algorithm is used with relative tolerance equal to
$10^{-6}$. The inhomogeneity $m(x)$ is confined in the unit cube and
the velocity field $c(x)$ for each example is between $2/3$ and
$1$. The grid step size $h$ is chosen such that there are $6$ points
per wavelength for the homogeneous region. This ensures a minimum of 4
points per wavelength across the whole domain.

\begin{figure}[h!]
  \begin{center}
    \begin{tabular}{|cc|cc|cc|cc|}
      \hline
      $\omega$ & $N$ & $T_s$(sec) & $T_a$(sec) & $n_p$ & $T_p$(sec) \\
      \hline
      1.0e+02 & 9.0e+03 & 1.3e+00 & 6.6e-02 & 5 & 4.2e-01\\
      2.0e+02 & 3.6e+04 & 6.4e+00 & 2.6e-01 & 5 & 1.4e+00\\
      4.0e+02 & 1.5e+05 & 3.9e+01 & 1.0e+00 & 6 & 7.8e+00\\
      \hline
    \end{tabular}
    \includegraphics[height=1.8in]{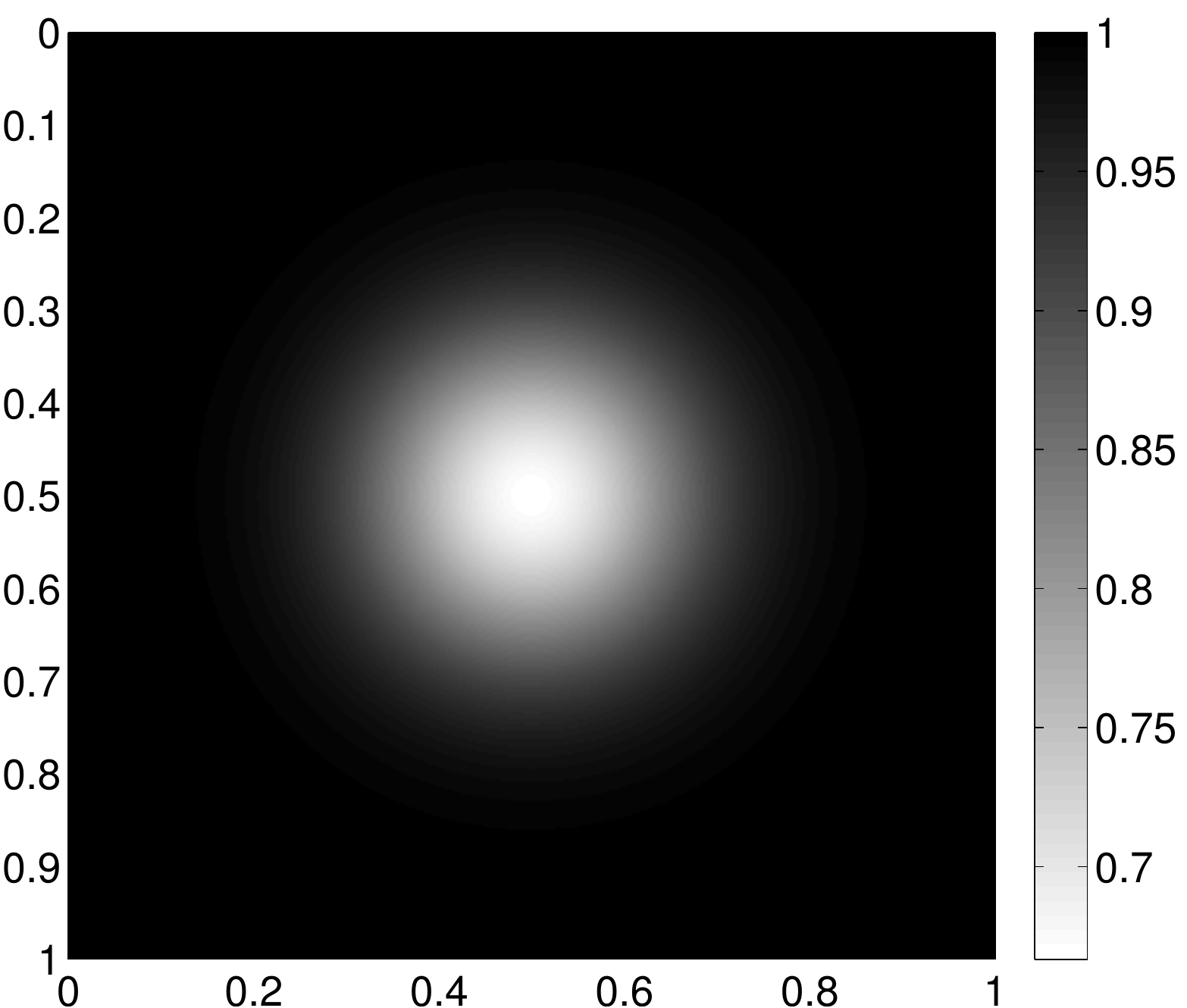} \hspace{0.25in} \includegraphics[height=1.8in]{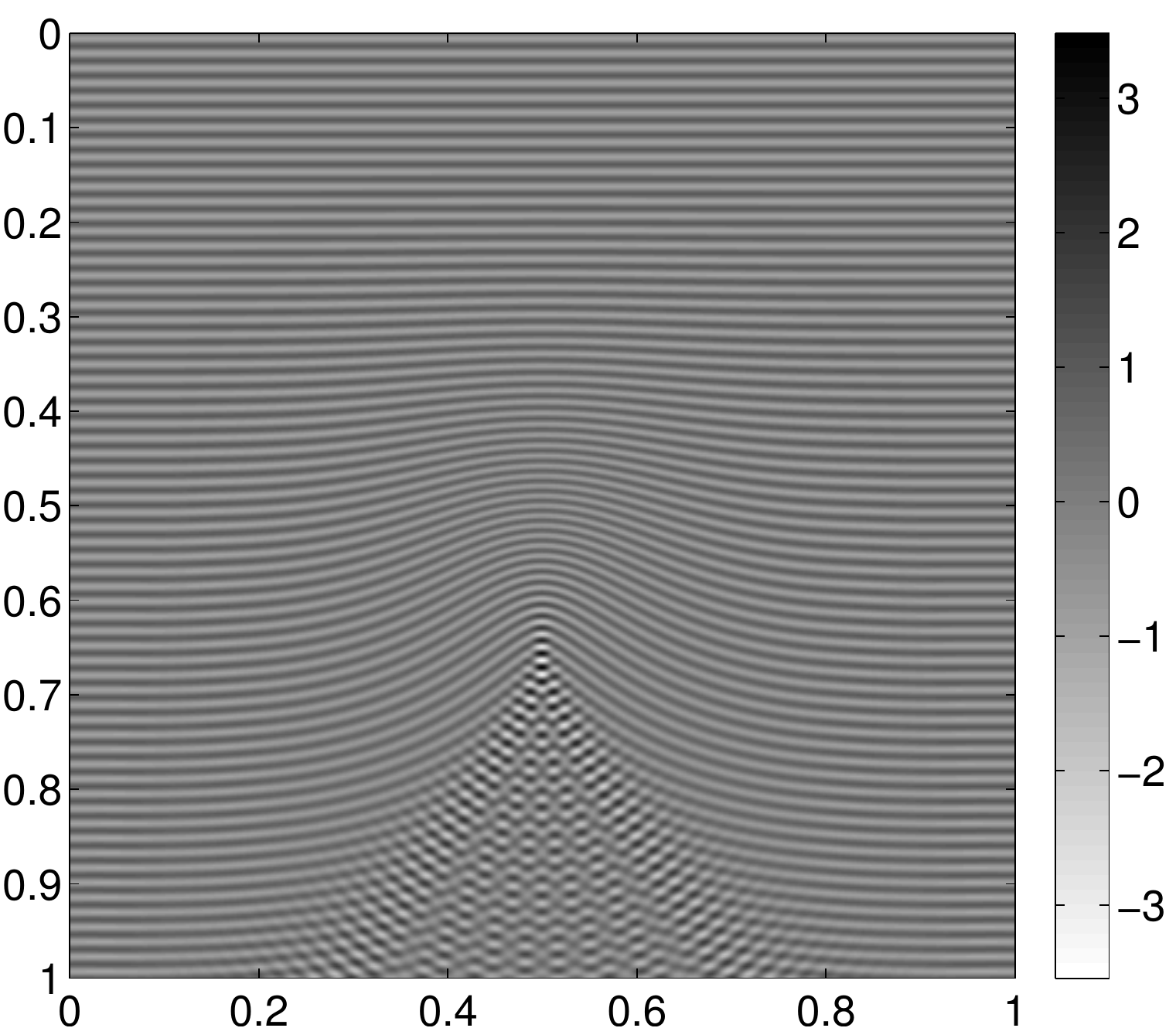}\\
  \end{center}
  \caption{Example 1 of 2D rectangular case. Top: numerical results.
    Bottom: $c(x)$ (left) and $u(x)+u_I(x)$ (right) for the largest
    $\omega$ value.}
  \label{fig:R21}
\end{figure}

\begin{figure}[h!]
  \begin{center}
    \begin{tabular}{|cc|cc|cc|cc|}
      \hline
      $\omega$ & $N$ & $T_s$(sec) & $T_a$(sec) & $n_p$ & $T_p$(sec) \\
      \hline
      1.0e+02 & 9.0e+03 & 1.3e+00 & 7.7e-02 & 6 & 4.3e-01\\
      2.0e+02 & 3.6e+04 & 7.0e+00 & 4.1e-01 & 8 & 4.1e+00\\
      4.0e+02 & 1.5e+05 & 4.2e+01 & 1.3e+00 & 8 & 1.3e+01\\
      \hline
    \end{tabular}
    \includegraphics[height=1.8in]{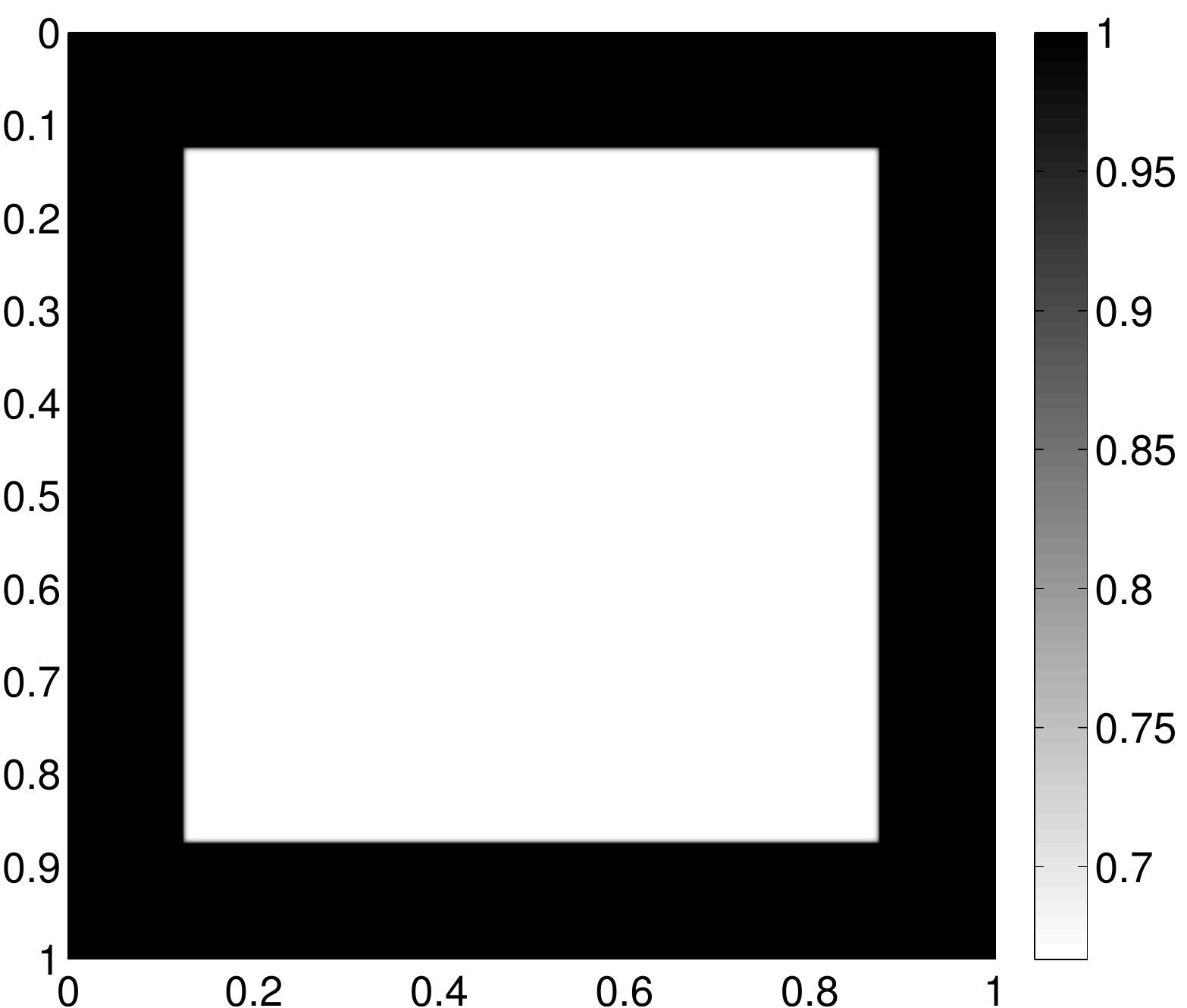} \hspace{0.25in} \includegraphics[height=1.8in]{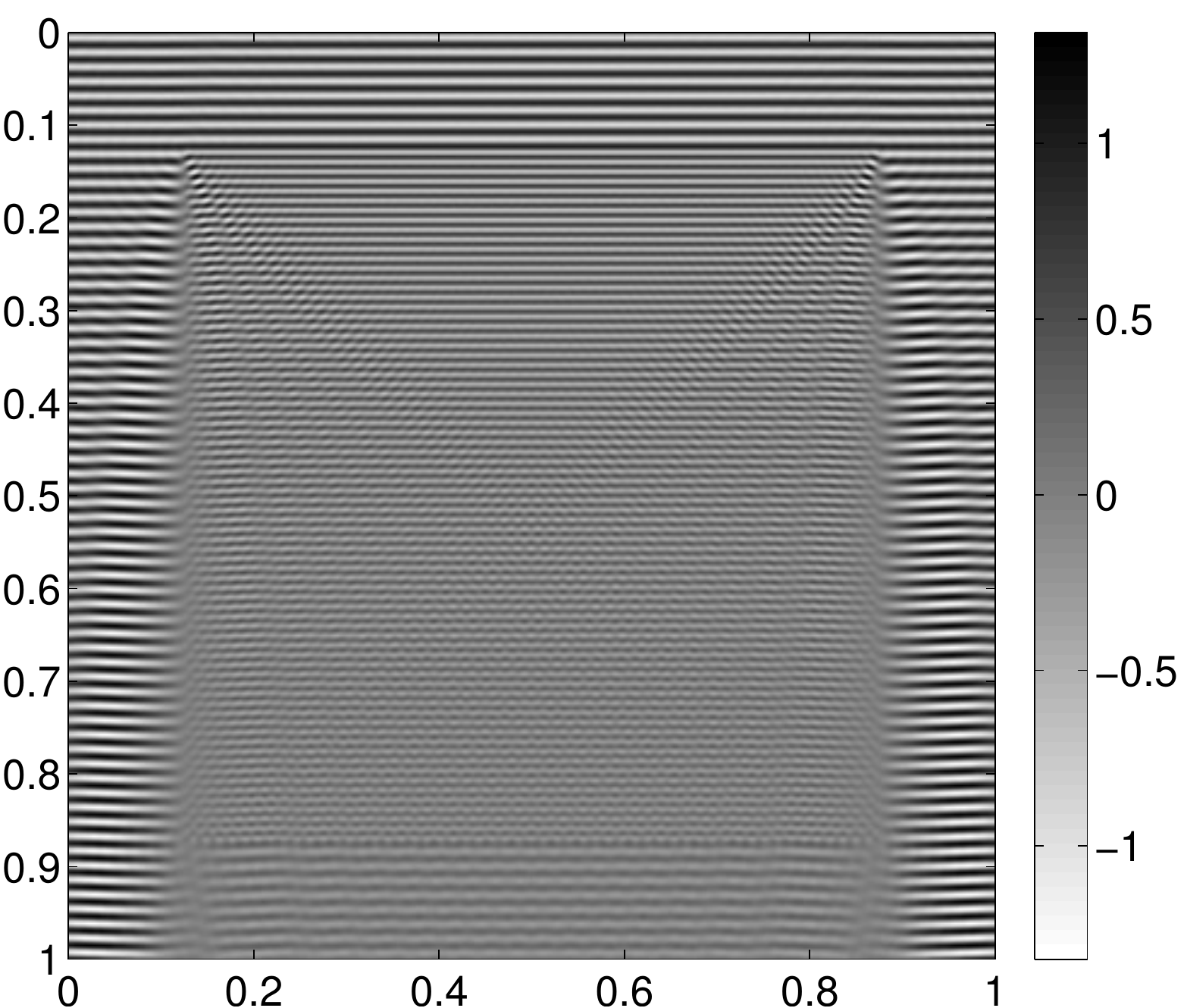}\\
  \end{center}
  \caption{Example 2 of 2D rectangular case. Top: numerical
    results. Bottom: $c(x)$ (left) and $u(x)+u_I(x)$ (right) for the
    largest $\omega$ value.}
  \label{fig:R22}
\end{figure}

The 2D case is tested with two examples: a Gaussian bump and a
smoothed square cavity. The incoming wave is a plane wave pointing
downward at frequency $\omega$ and the results are summarized in
Figures \ref{fig:R21} and \ref{fig:R22}. The columns of the tables
are:
\begin{itemize}
\item $\omega$ is the frequency,
\item $N$ is the number of unknowns,
\item $T_s$ is the setup time of the preconditioner in seconds,
\item $T_a$ is the application time of the preconditioner in seconds,
\item $n_p$ is the iteration number of the preconditioned iteration, and
\item $T_p$ is the solution time of the preconditioned iteration in
  seconds.
\end{itemize}

\begin{figure}[h!]
  \begin{center}
    \begin{tabular}{|cc|cc|cc|}
      \hline
      $\omega$ & $N$ & $T_s$(sec) & $T_a$(sec) & $n_p$ & $T_p$(sec) \\
      \hline
      2.5e+01 & 1.2e+04 & 5.2e+00 & 6.3e-02 & 4 & 2.8e-01\\
      5.0e+01 & 1.0e+05 & 1.0e+02 & 6.4e-01 & 4 & 2.9e+00\\
      1.0e+02 & 8.6e+05 & 3.2e+03 & 7.8e+00 & 4 & 3.4e+01\\
      \hline
    \end{tabular}
    \includegraphics[height=1.8in]{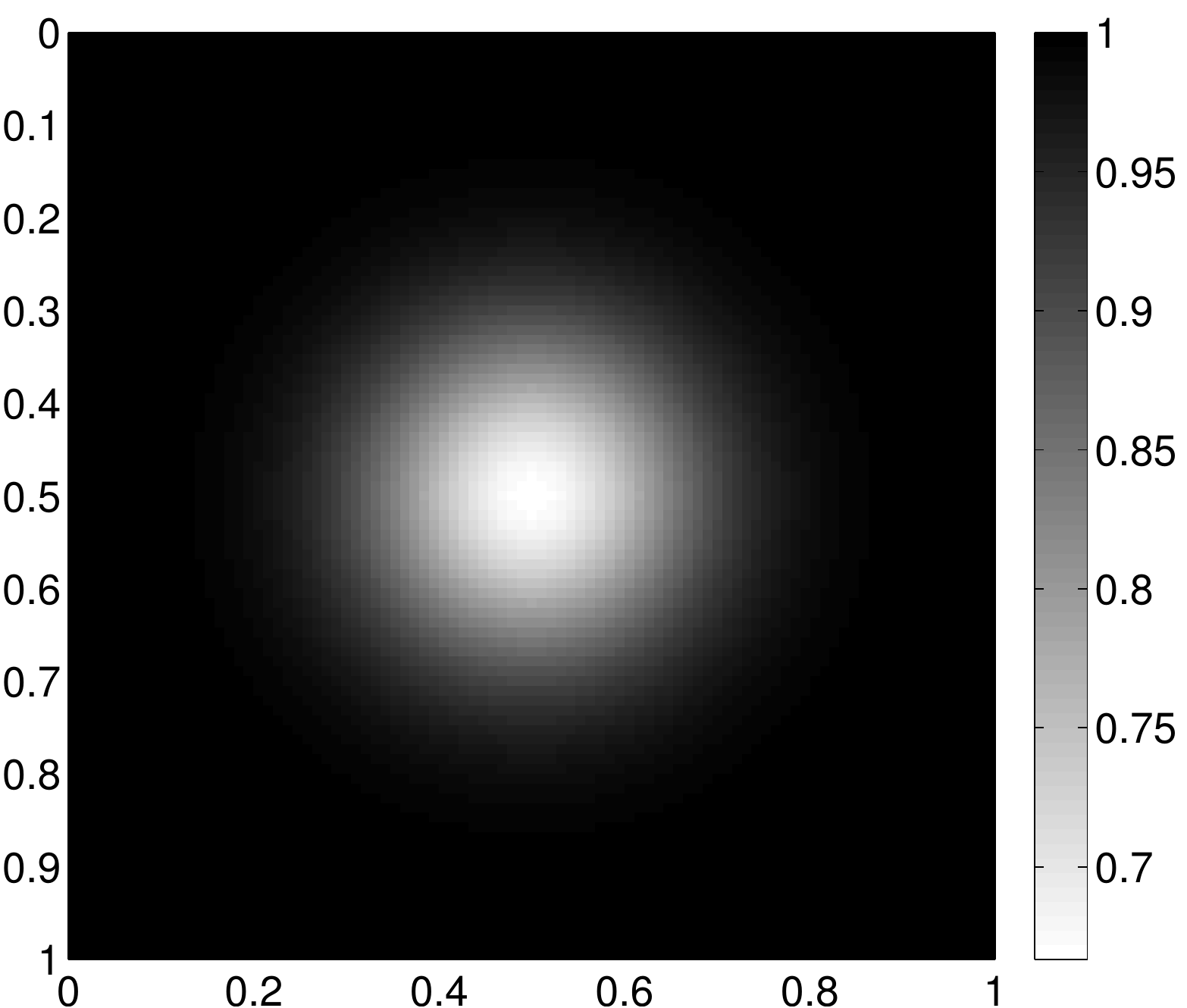} \hspace{0.25in} \includegraphics[height=1.8in]{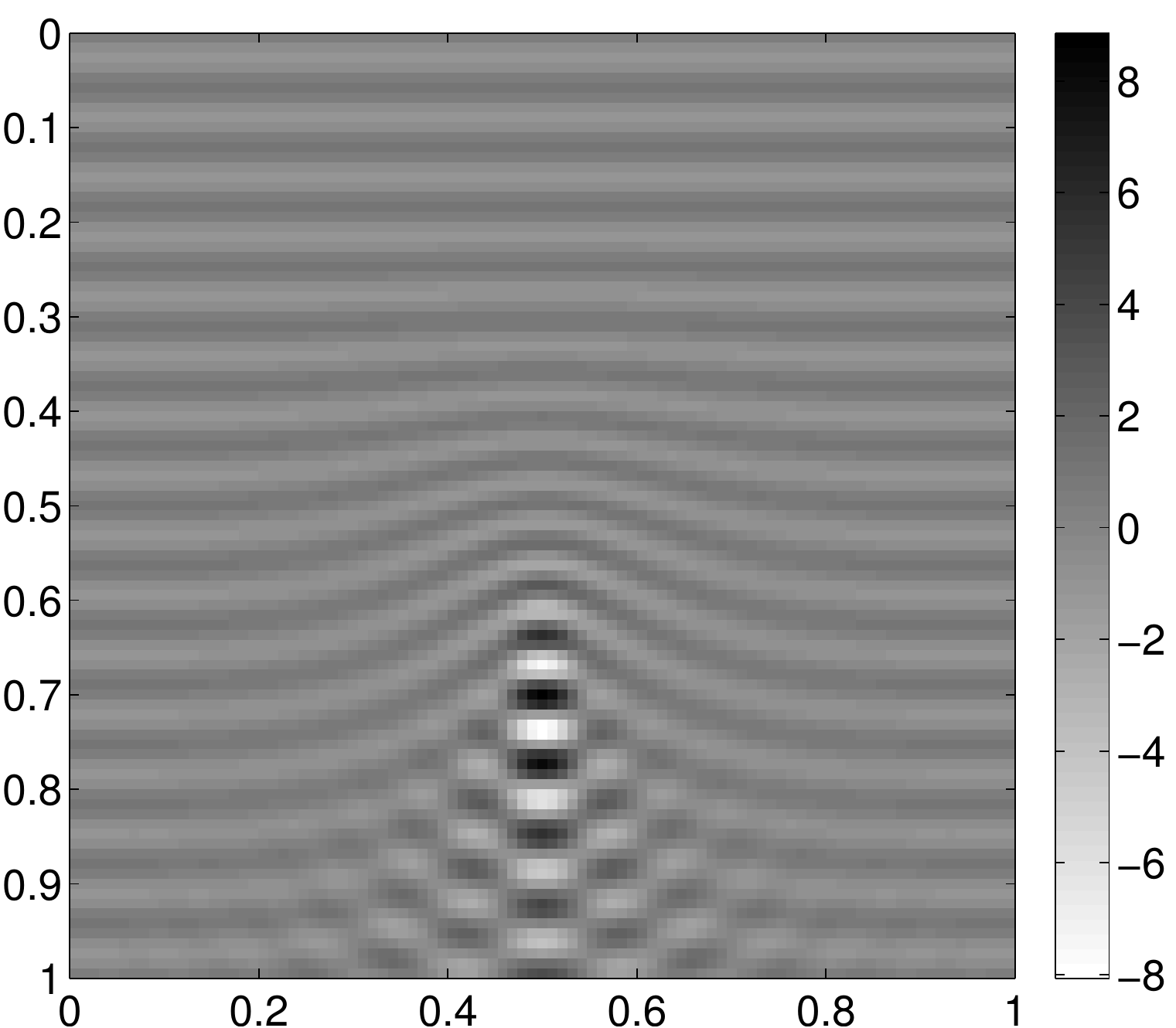}\\
  \end{center}
  \caption{Example 1 of 3D rectangular case. Top: numerical
    results. Bottom: cross-sections of $c(x)$ (left) and $u(x)+u_I(x)$
    (right) in the middle of the domain for the largest $\omega$
    value.}
  \label{fig:R31}
\end{figure}

\begin{figure}[h!]
  \begin{center}
    \begin{tabular}{|cc|cc|cc|}
      \hline
      $\omega$ & $N$ & $T_s$(sec) & $T_a$(sec) & $n_p$ & $T_p$(sec) \\
      \hline
      2.5e+01 & 1.2e+04 & 5.6e+00 & 7.3e-02 & 6 & 4.7e-01\\
      5.0e+01 & 1.0e+05 & 1.0e+02 & 5.9e-01 & 8 & 5.7e+00\\
      1.0e+02 & 8.6e+05 & 2.9e+03 & 7.7e+00 & 9 & 7.5e+01\\
      \hline
    \end{tabular}
    \includegraphics[height=1.8in]{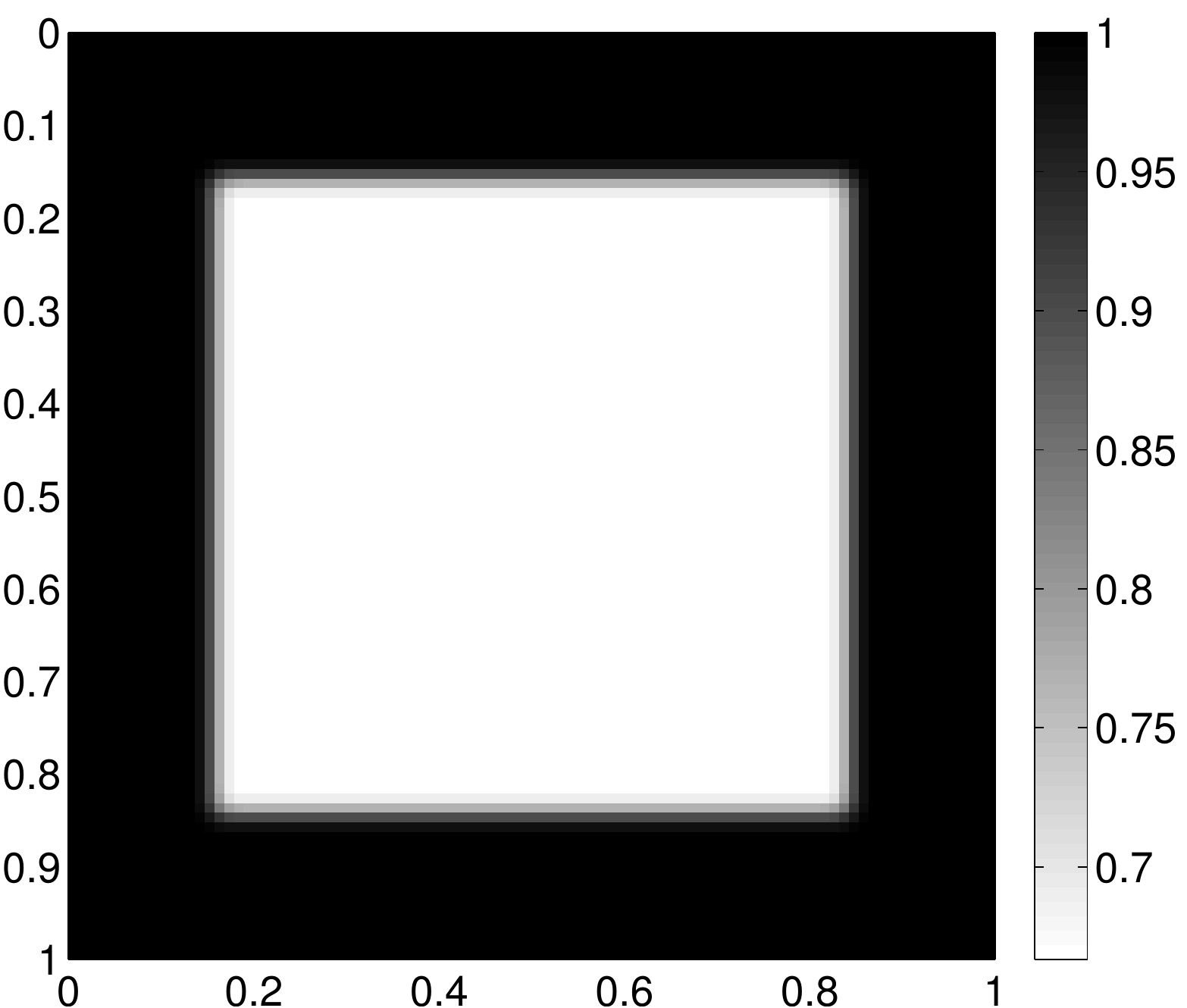} \hspace{0.25in} \includegraphics[height=1.8in]{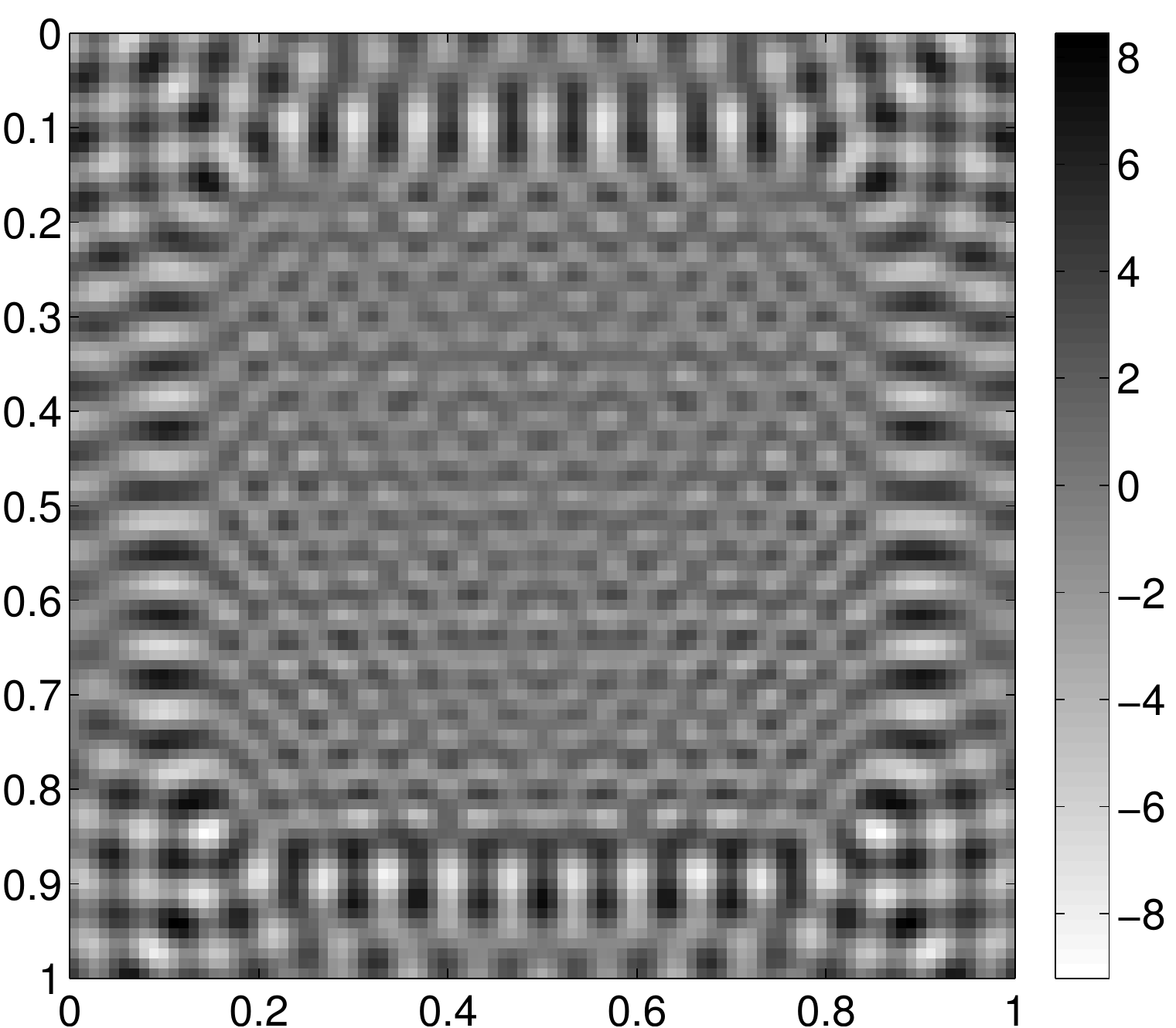}\\
  \end{center}i
  \caption{Example 2 of 3D rectangular case. Top: numerical
    results. Bottom: cross-sections of $c(x)$ (left) and $u(x)+u_I(x)$
    (right) in the middle of the domain for the largest $\omega$
    value.}
  \label{fig:R32}
\end{figure}

The 3D case is also tested with two examples: a Gaussian bump and a
smoothed cubic cavity. The incoming wave is again a plane wave
pointing downward at frequency $\omega$ and the results are given in
Figures \ref{fig:R31} and \ref{fig:R32}.

The results show that the setup and application costs of the
preconditioner scale with $\omega$ and $N$ according to the complexity
analysis given above. The preconditioner reduces the iteration number
dramatically and in fact it becomes essentially frequency-independent.

\section{General domains}
\label{sec:G}

For most problems, the support $\Omega$ of $m(x)$ is not necessarily
rectangular. We can apply the approach in Section \ref{sec:R} by
embedding $\Omega$ into a rectangular domain. However, this approach
increases the number of unknowns significantly, especially in 3D. This
section proposes a different approach that does not suffer from this.

\subsection{Algorithm}

The construction of the operator $A$ for the interior grid points is
exactly the same as the one in Section \ref{sec:R}. The main
difference is in the construction of $B$. For a boundary point $i$, we
require
\[
B(i,\mu(i)) (K\omega^2 m)(\mu(i),:) \approx 0.
\]
Since the local geometry at $i$ can be quite different from point to
point, it is impossible to find a uniform stencil $B(i,\mu(i))$ as was
done in Section \ref{sec:R}. Therefore, one needs to consider $i$
point by point. However, applying the approach of Section \ref{sec:R}
to each $i$ individually is too costly since each singular value
decomposition costs $O(N)$ steps. To keep the complexity under
control, we propose the following randomized approach.

Let $R$ be a Gaussian random matrix of size $N\times r$ where $r=O(1)$
is a constant multiple of the maximum stencil size (i.e., $9$ in 2D
and $27$ in 3D) and define the $N\times r$ matrix $T$ via
\[
T = (K \omega^2 m) R.
\]
For each boundary $i$, we look for $B(i,\mu(i))$ such that
\[
B(i,\mu(i))  T(\mu(i),:) \approx 0.
\]
This is done by solving the optimization problem
\[
\min_{\alpha: \|\alpha\|=1} \|\alpha\cdot T(\mu(i), :)\|^2.
\]
With the singular value decomposition $T(\mu(i),:) = USV^*$, the optimal $\alpha$ is given by
\[
\alpha = U(:,|\mu(i)|)^*.
\]
Finally, we set $B(i,\mu(i)) = \alpha$. This process is repeated for
each boundary point $i$. Since $T(\mu(i),:)$ is of size $O(1)\times
O(1)$, the computational cost of the singular value decomposition setup for each $i$ is also of
order $O(1)$.

The effectiveness of this preconditioner depends on the assumption
that there exists a good local stencil in $\mu(i)$ for a boundary
point $i$. When the domain is convex or nearly convex, such a stencil
is guaranteed due to the absorbing-type boundary conditions.  However,
when the domain is highly non-convex the preconditioner is less
effective.

\subsection{Complexity}
To analyze the cost of this approach, we assume that the size of the
support of $m(x)$ is more than a constant fraction of the size of its
bounding box, which implies that $N=O(n^d)$.

In 2D, the setup algorithm consists of three parts: (i) forming $T =
(K \omega^2 m) R$, (ii) computing an singular value decomposition for
each boundary point $i$, and (iii) factorizing \eqref{eq:ABop} with
the nested dissection algorithm.
\begin{itemize}
\item The first part can be accomplished in $O(N\log N)$ steps by
  using the fast Fourier transform.
\item The second part takes $O(\sqrt{N})$ steps since there are at
  most $O(\sqrt{N})$ boundary points and computing singular value
  decomposition for each one takes $O(1)$ steps.
\item Finally, factorizing with the nested dissection algorithm takes
  $O(N^{3/2})$ steps.
\end{itemize}
Adding these together shows that the construction cost is
$O(N^{3/2})$. The application cost of the preconditioner is $O(N\log
N)$ since it is essentially a solve with an existing 2D nested
dissection factorization.

In 3D, among the three parts of the setup algorithm, only the cost of
computing a nested dissection factorization increases to $O(N^2)$. This
now becomes the dominant part of the setup cost. The application
algorithm is again a solve with the existing nested dissection
factorization with a cost of $O(N^{4/3})$. 

Though the asymptotic complexity obtained for this case is similar to
the one for the rectangular case, the actual running time is often
lower because the sizes of the supernodes in the nested dissection
algorithm is often much smaller.

\subsection{Numerical results}

\begin{figure}[h!]
  \begin{center}
    \begin{tabular}{|cc|cc|cc|}
      \hline
      $\omega$ & $N$ & $T_s$(sec) & $T_a$(sec) & $n_p$ & $T_p$(sec) \\
      \hline
      1.0e+02 & 7.6e+03 & 9.8e-01 & 4.1e-02 & 6 & 2.4e-01\\
      2.0e+02 & 3.0e+04 & 3.7e+00 & 1.3e-01 & 7 & 1.2e+00\\
      4.0e+02 & 1.2e+05 & 2.1e+01 & 5.1e-01 & 8 & 6.5e+00\\
      \hline
    \end{tabular}
    \includegraphics[height=1.8in]{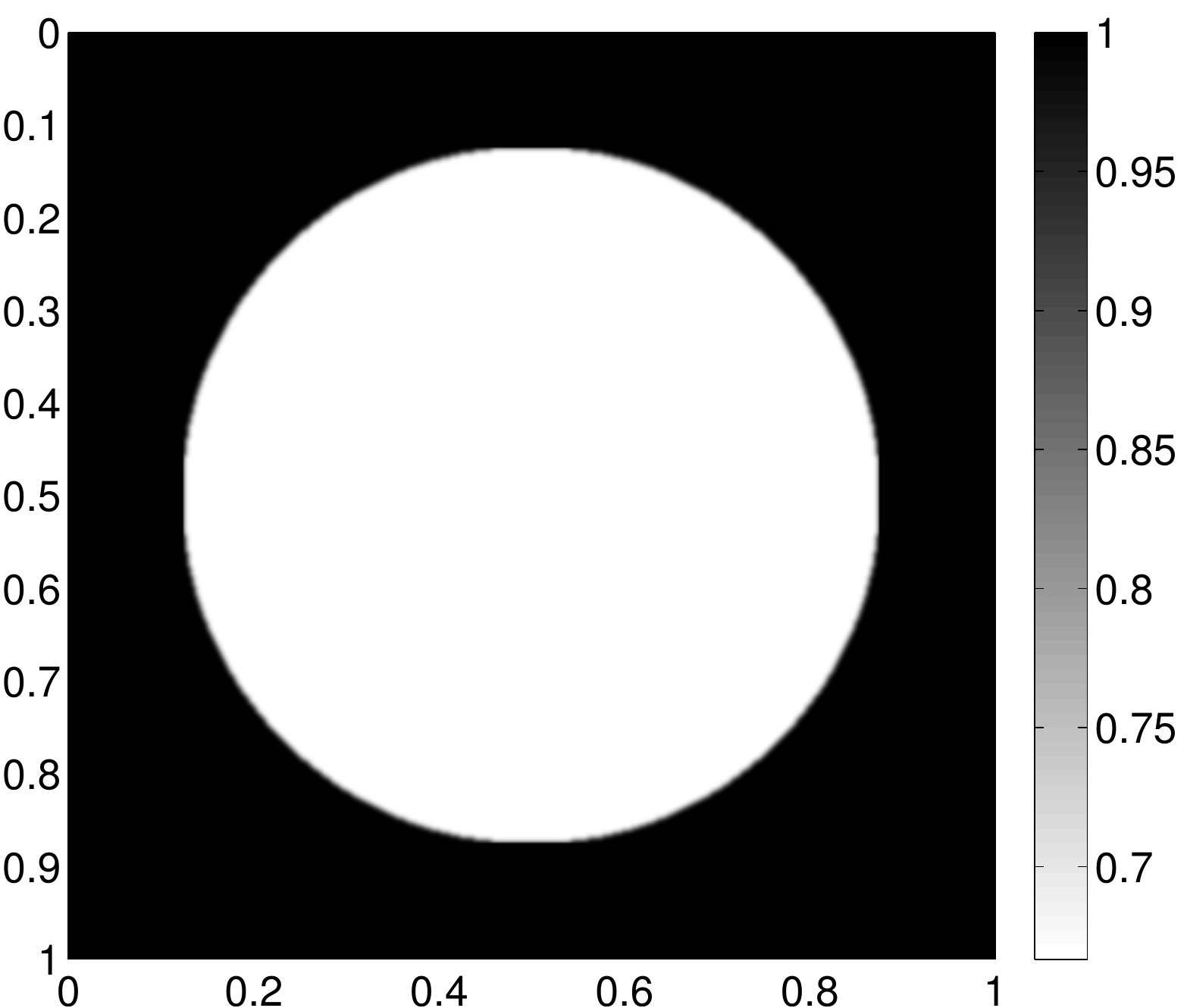} \hspace{0.25in} \includegraphics[height=1.8in]{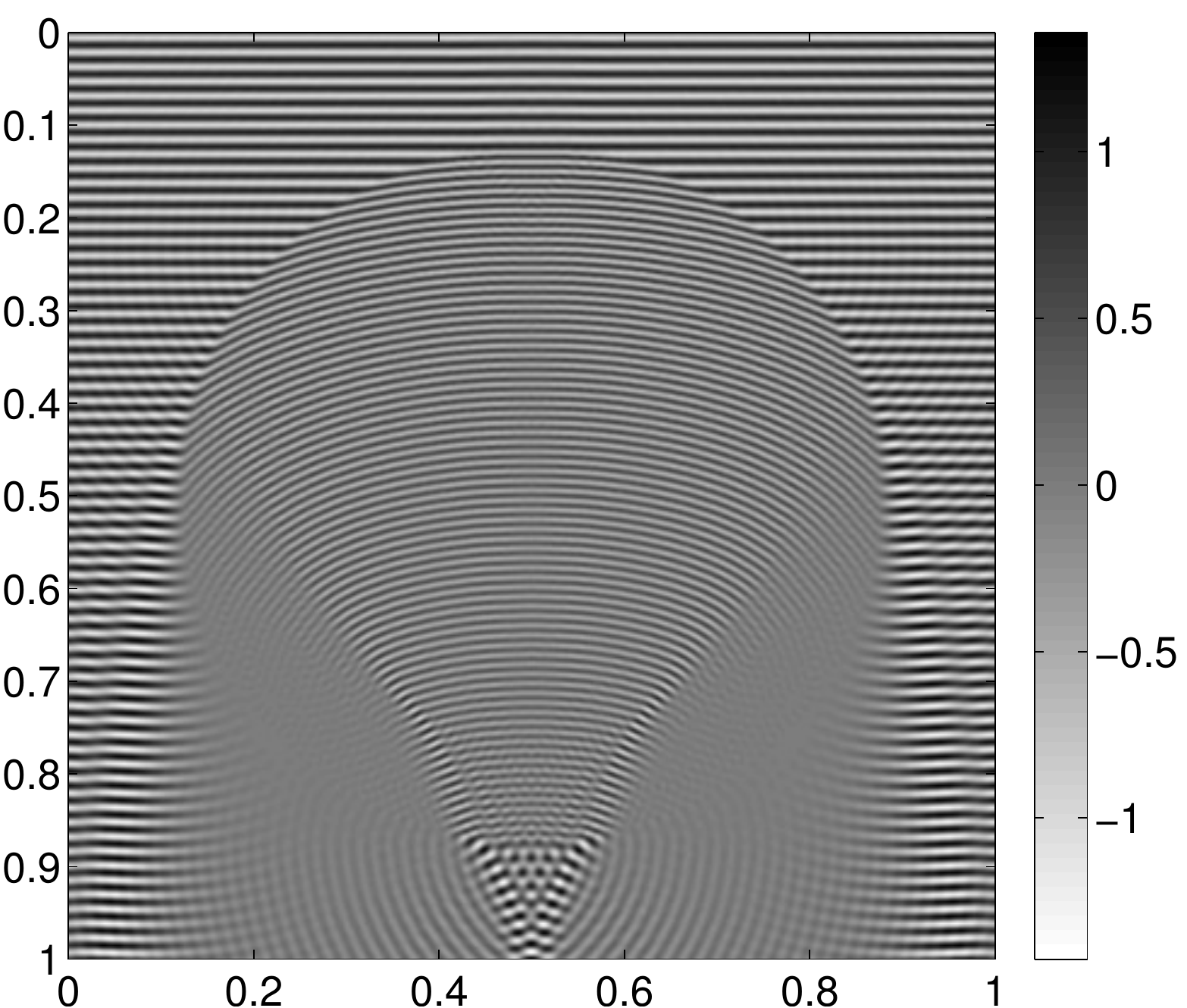}\\
  \end{center}
  \caption{Example 1 of 2D general case. Top: numerical results.
    Bottom: $c(x)$ (left) and $u(x)+u_I(x)$ (right) for the largest
    $\omega$ value.}
  \label{fig:G21}
\end{figure}

\begin{figure}[h!]
  \begin{center}
    \begin{tabular}{|cc|cc|cc|}
      \hline
      $\omega$ & $N$ & $T_s$(sec) & $T_a$(sec) & $n_p$ & $T_p$(sec) \\
      \hline
      1.0e+02 & 5.7e+03 & 7.4e-01 & 3.3e-02 & 5 & 2.0e-01\\
      2.0e+02 & 2.1e+04 & 2.8e+00 & 7.6e-02 & 8 & 1.1e+00\\
      4.0e+02 & 7.9e+04 & 1.8e+01 & 2.9e-01 & 8 & 5.3e+00\\
      \hline
    \end{tabular}
    \includegraphics[height=1.8in]{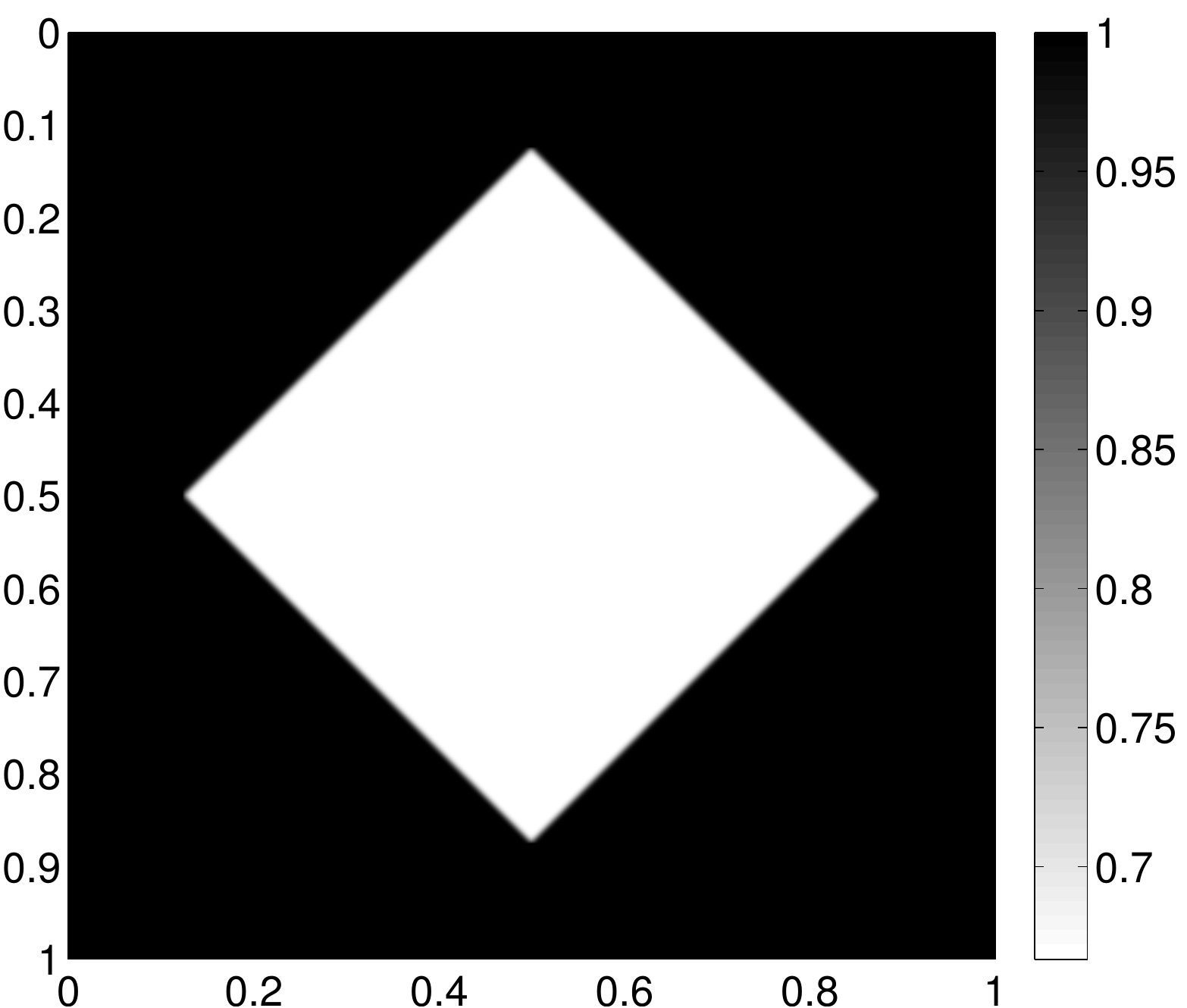} \hspace{0.25in} \includegraphics[height=1.8in]{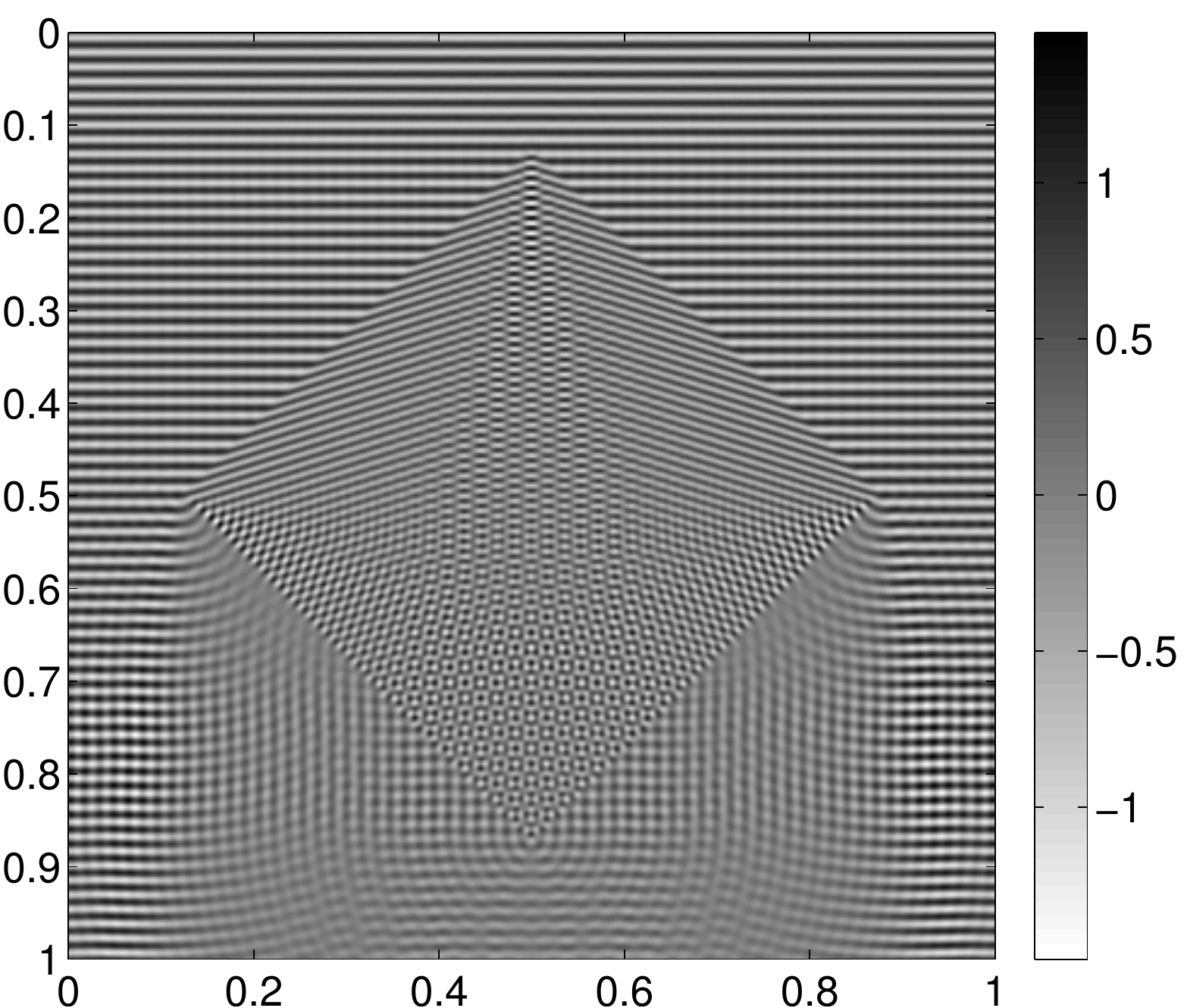}\\
  \end{center}
  \caption{Example 2 of 2D general case. Top: numerical results.
    Bottom: $c(x)$ (left) and $u(x)+u_I(x)$ (right) for the largest
    $\omega$ value.}
  \label{fig:G22}
\end{figure}

The 2D case is tested with two examples: a smoothed cavity of an
$\ell_2$ ball and a smoothed cavity of an $\ell_1$ ball. The incoming
wave is a plane wave pointing downward at frequency $\omega$ and the
numerical results are given in Figures \ref{fig:G21} and
\ref{fig:G22}.

\begin{figure}[h!]
  \begin{center}
    \begin{tabular}{|cc|cc|cc|}
      \hline
      $\omega$ & $N$ & $T_s$(sec) & $T_a$(sec) & $n_p$ & $T_p$(sec) \\
      \hline
      2.5e+01 & 6.5e+03 & 7.1e+00 & 2.7e-02 & 4 & 1.6e-01\\
      5.0e+01 & 5.5e+04 & 7.9e+01 & 3.4e-01 & 5 & 2.1e+00\\
      1.0e+02 & 4.5e+05 & 1.4e+03 & 4.1e+00 & 5 & 2.3e+01\\
      \hline
    \end{tabular}
    \includegraphics[height=1.8in]{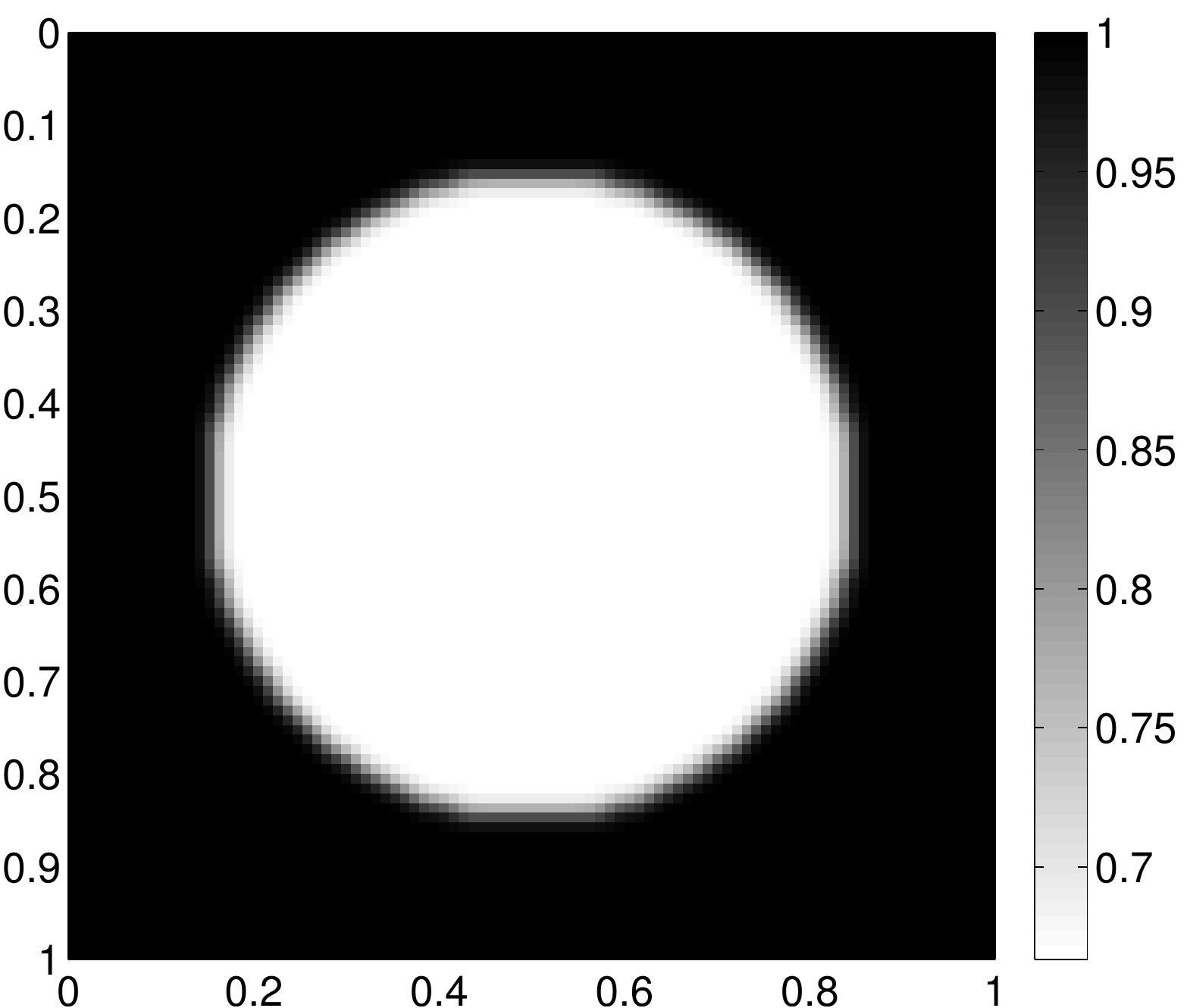} \hspace{0.25in} \includegraphics[height=1.8in]{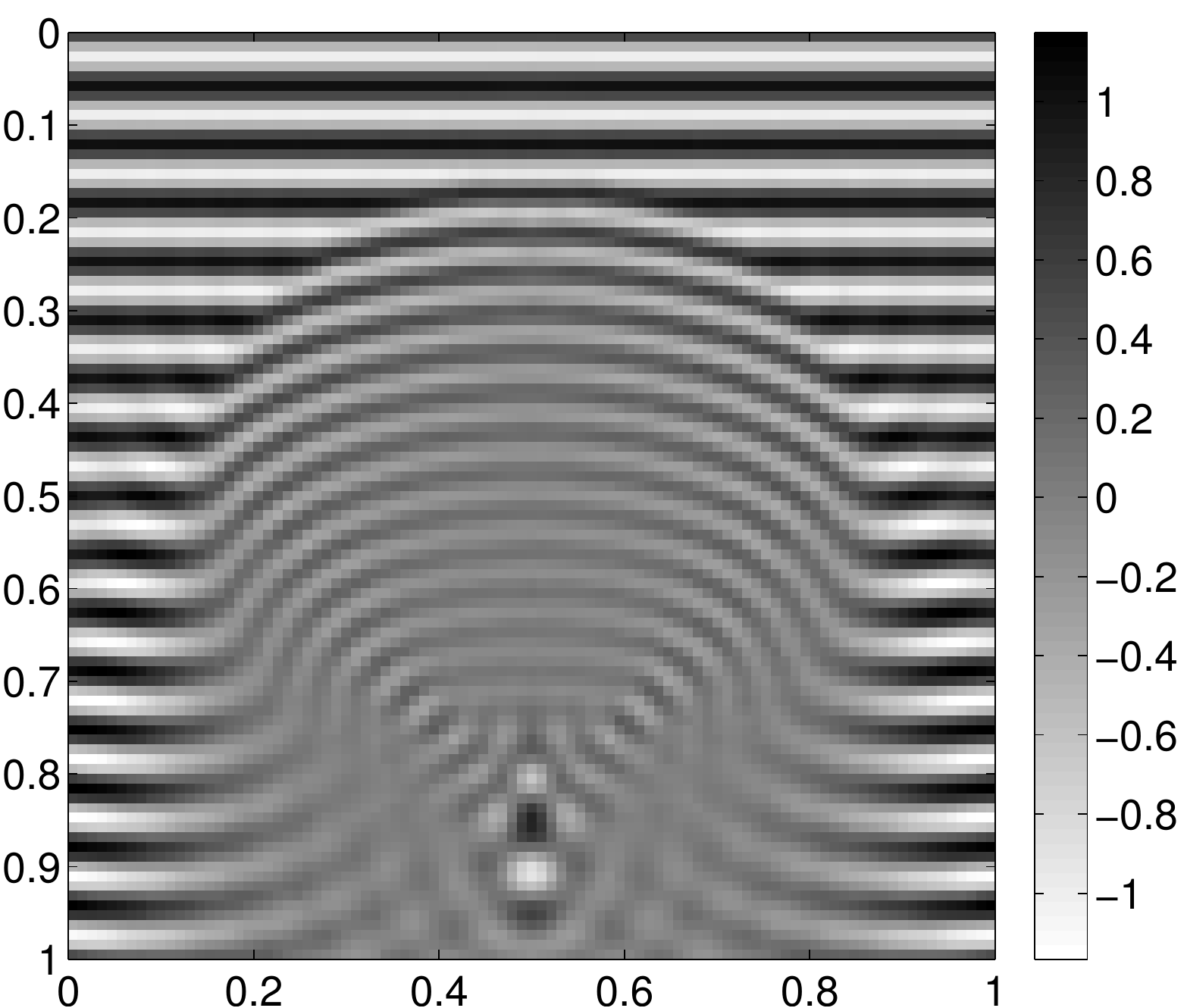}\\
  \end{center}
  \caption{Example 1 of 3D general case. Top: numerical
    results. Bottom: cross-sections of $c(x)$ (left) and $u(x)+u_I(x)$
    (right) in the middle of the domain for the largest $\omega$
    value.}
  \label{fig:G31}
\end{figure}

\begin{figure}[h!]
  \begin{center}
    \begin{tabular}{|cc|cc|cc|}
      \hline
      $\omega$ & $N$ & $T_s$(sec) & $T_a$(sec) & $n_p$ & $T_p$(sec) \\
      \hline
      2.5e+01 & 3.7e+03 & 4.9e+00 & 3.0e-02 & 4 & 1.0e-01\\
      5.0e+01 & 3.1e+04 & 4.5e+01 & 1.8e-01 & 6 & 1.4e+00\\
      1.0e+02 & 2.3e+05 & 7.5e+02 & 2.1e+00 & 9 & 2.3e+01\\
      \hline
    \end{tabular}
    \includegraphics[height=1.8in]{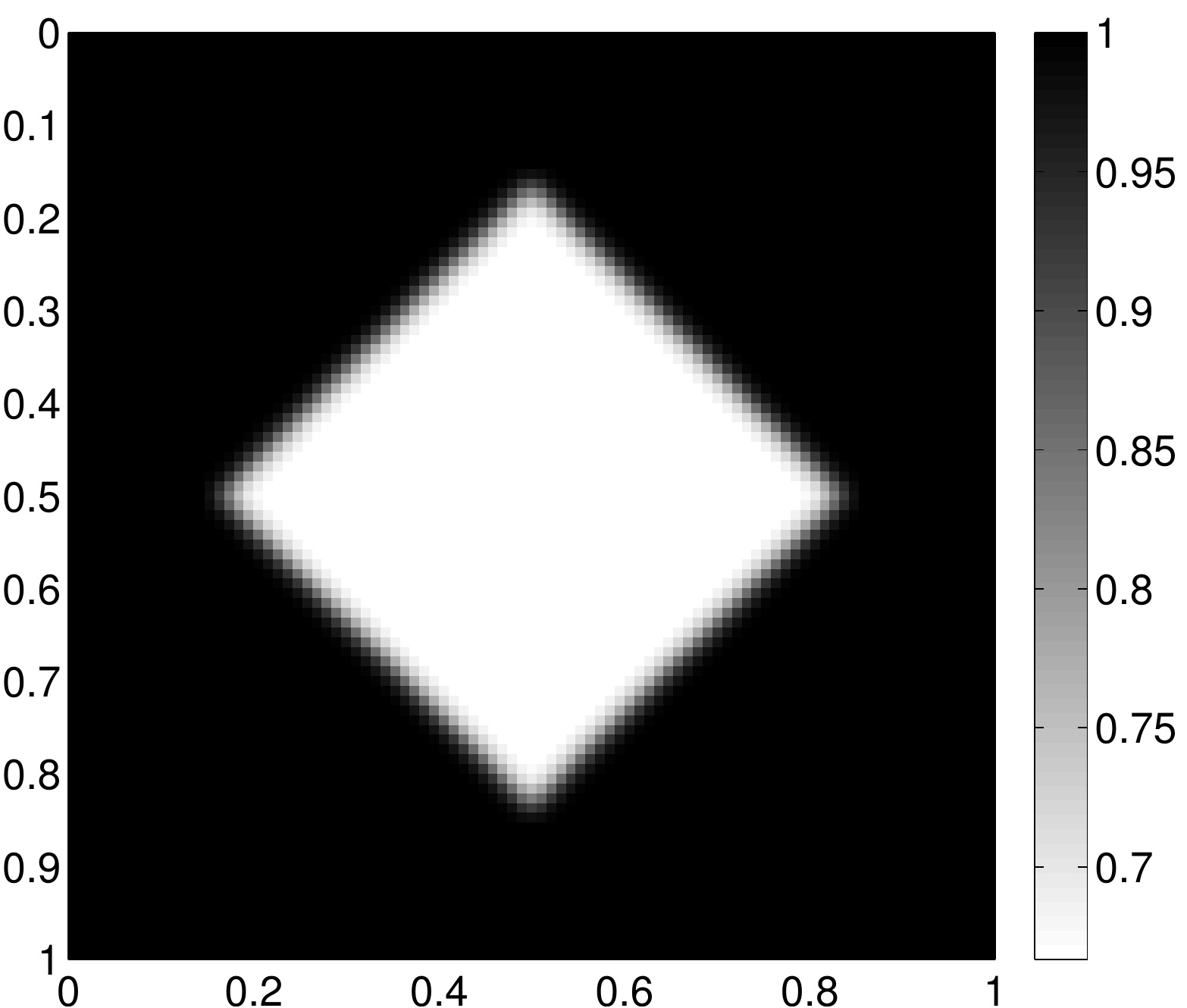} \hspace{0.25in} \includegraphics[height=1.8in]{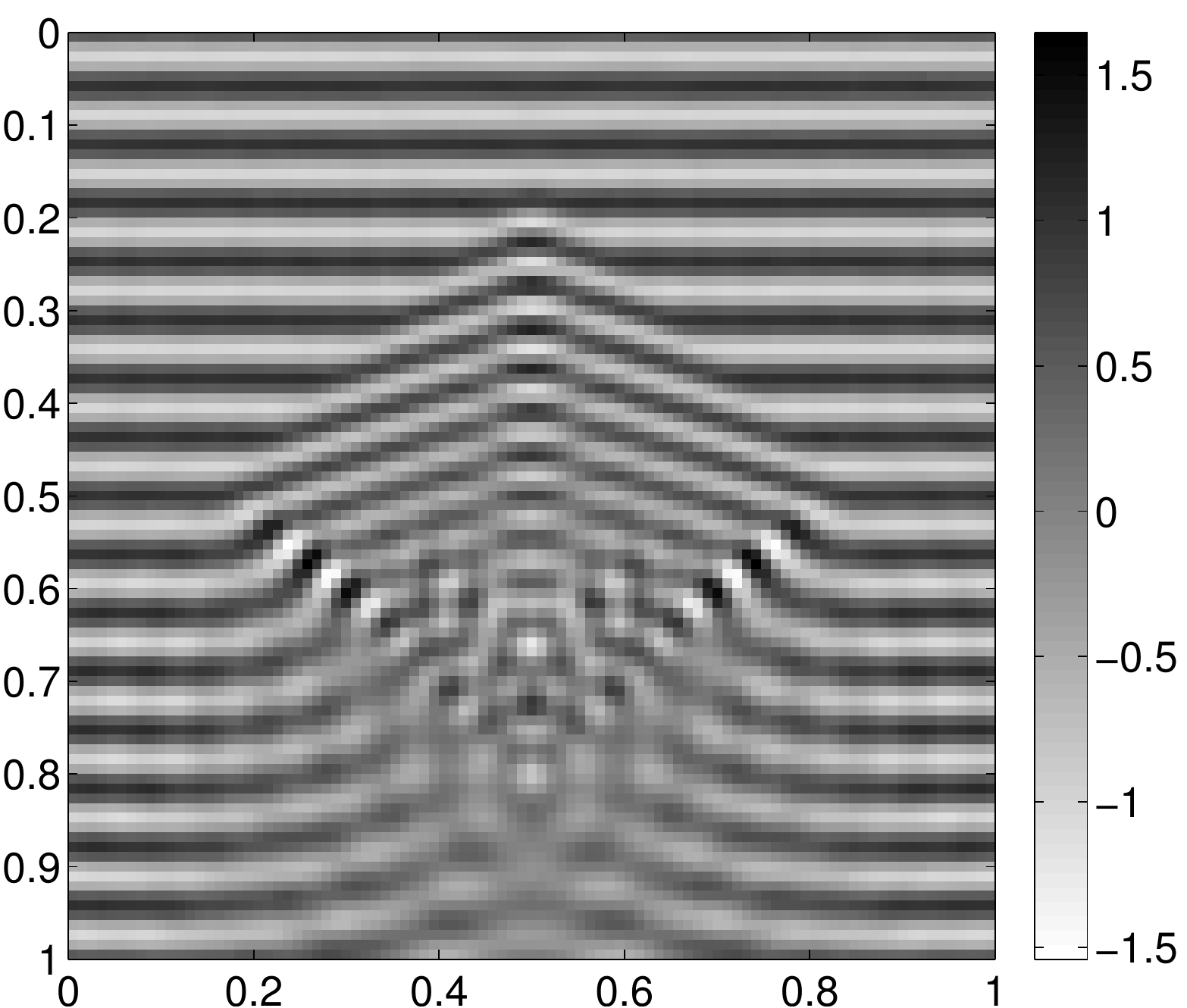}\\
  \end{center}
  \caption{Example 2 of 3D general case. Top: numerical
    results. Bottom: cross-sections of $c(x)$ (left) and $u(x)+u_I(x)$
    (right) in the middle of the domain for the largest $\omega$
    value.}
  \label{fig:G32}
\end{figure}

Two examples are tested for the 3D case: a smoothed cavity of an
$\ell_2$ ball and a smoothed cavity of an $\ell_1$ ball. The incoming
wave is a plane wave pointing downward at frequency $\omega$ and the
numerical results are given in Figures \ref{fig:G31} and
\ref{fig:G32}.

The results show that the setup and application costs of the
preconditioner are lower compared to the ones of the rectangular
case. The number of iterations remains essentially independent of the
frequency, indicating that the boundary stencils computed randomly
provide a sufficiently accurate boundary condition for the solution.

\section{Laplace equation}

In this section, we consider the 3D Laplace equation with compact
potential perturbation, i.e.,
\begin{equation}
  (-\Delta + V(x))u(x) = f(x),\quad x\in\R^3
  \label{eq:lap}
\end{equation}
where $V(x)$ is supported in $\Omega = (0,1)^3$ and $u(x)$ decays like
$1/|x|$ as $x$ goes to infinity. This can be regarded as a limiting
case as $\omega$ goes to zero. The Green's function is given by
\[
G(x) = \frac{1}{4 \pi |x|}.
\]
Convolving $G(x)$ with \eqref{eq:lap} gives the Lippmann-Schwinger
equation
\begin{equation}
  u + G\ast(Vu) = G\ast f.
  \label{eq:lslap}
\end{equation}
For a potential $V$ that is not too negative, this problem remains
elliptic and standard iteration methods such as GMRES converge within
a small number of iterations. However, when $V(x)$ is sufficiently
negative, the problem loses ellipticity and the iteration number
increases dramatically. 

As we mentioned, when $\omega$ is finite, PML and ABC offer reasonable
approximations to the Sommerfeld radiation condition. When $\omega$
goes to zero, it is not clear how to define analytically a similar
local boundary condition at a finite distance for the $1/|x|$ decaying
condition. However, the sparsifying preconditioner proposed in
Sections \ref{sec:R} and \ref{sec:G} can be applied to
\eqref{eq:lslap} without any modification and the resulting operator
$B$ offers a good numerical approximation to the $1/|x|$ decaying
condition.

\subsection{Numerical results}
\begin{figure}[h!]
  \begin{center}
    \begin{tabular}{|cc|cc|cc|}
      \hline
      $\max(|V|)$ & $N$ & $T_s$(sec) & $T_a$(sec) & $n_p$ & $T_p$(sec) \\
      \hline
      5.8e+02 & 1.2e+04 & 5.0e+00 & 4.9e-02 & 5 & 2.1e-01\\
      2.3e+03 & 1.0e+05 & 9.8e+01 & 2.8e-01 & 5 & 1.8e+00\\
      9.2e+03 & 8.6e+05 & 2.1e+03 & 3.0e+00 & 6 & 2.1e+01\\
      \hline
    \end{tabular}
    \includegraphics[height=1.8in]{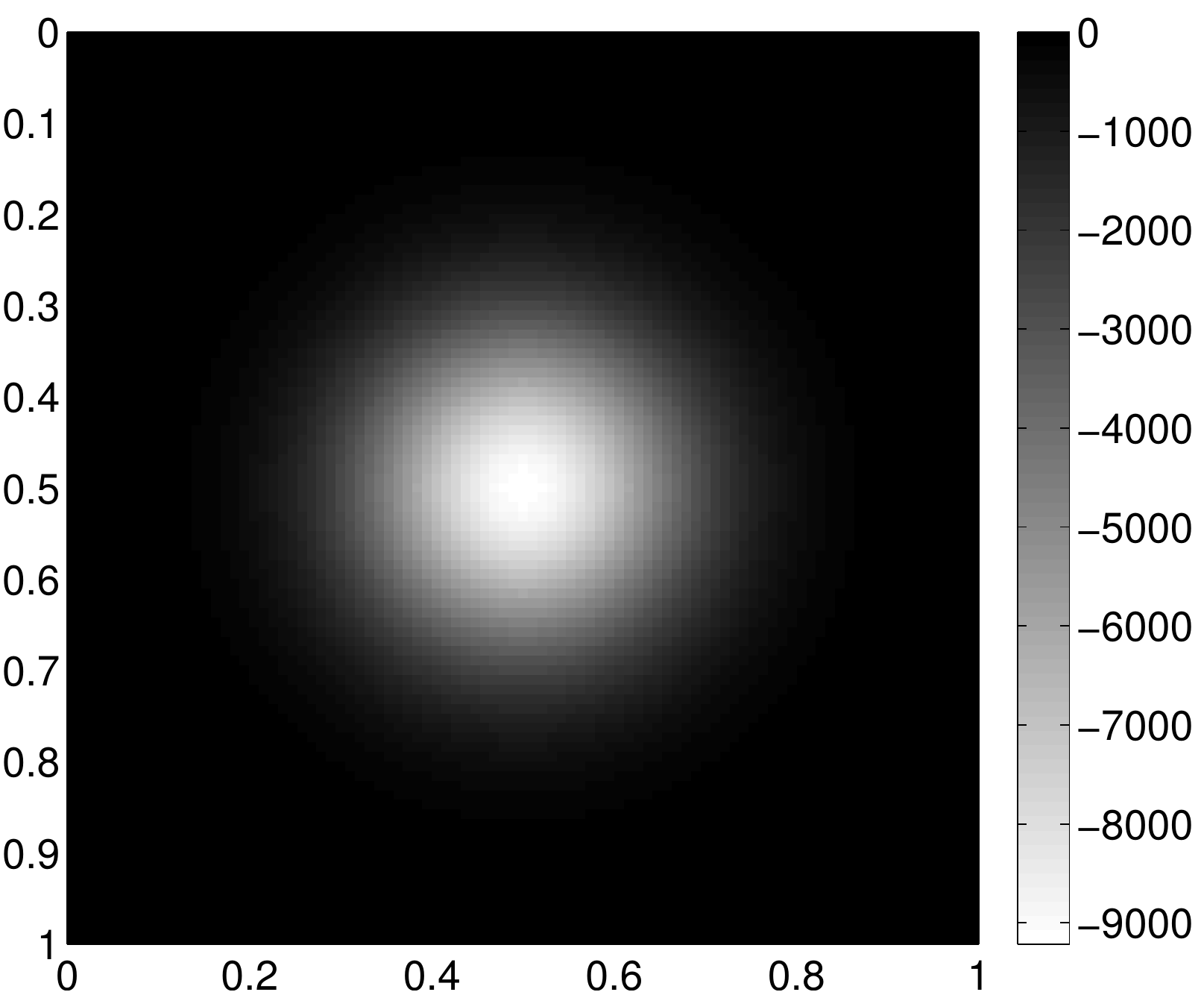} \hspace{0.25in} \includegraphics[height=1.8in]{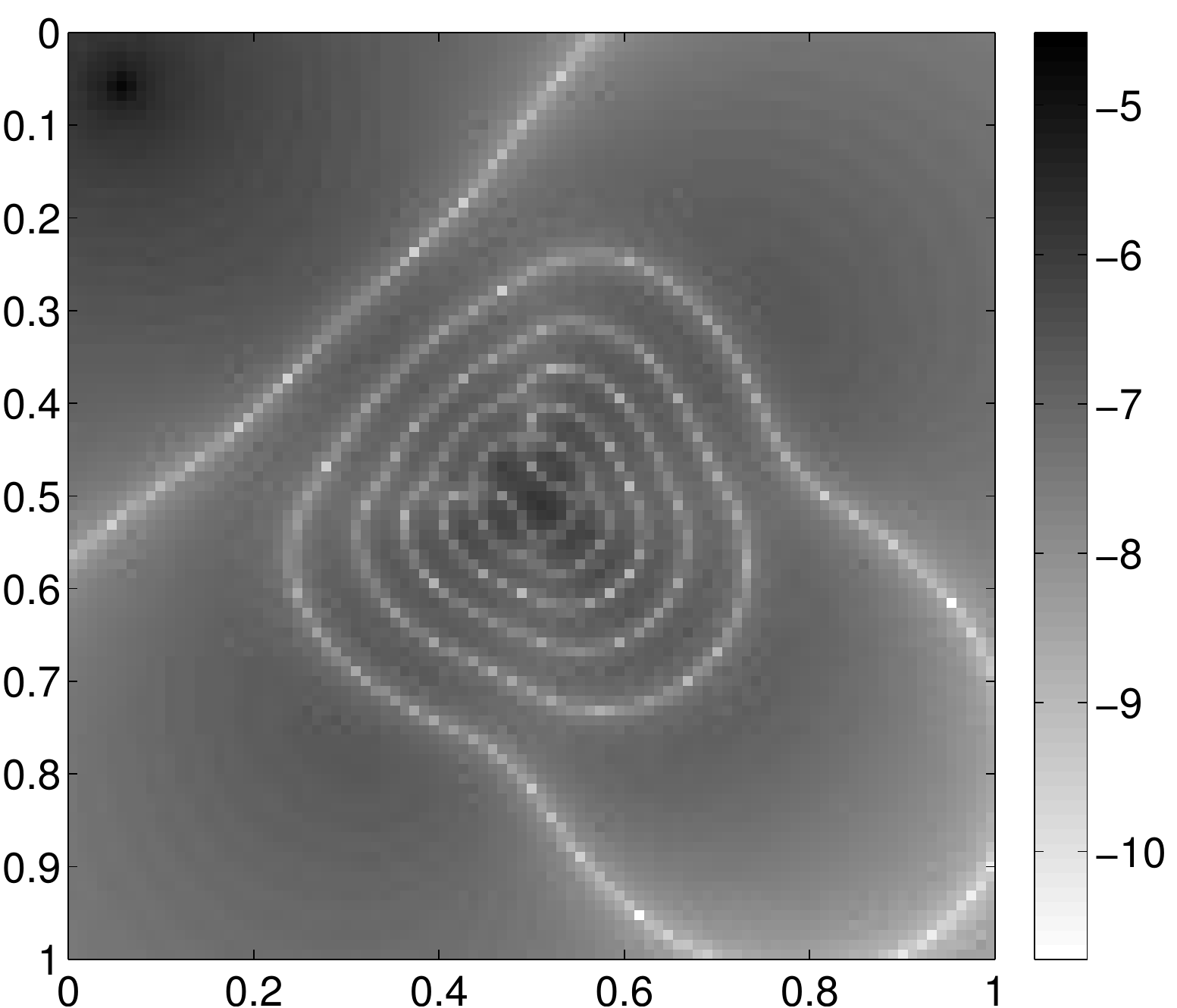}\\
  \end{center}
  \caption{Example 1 of 3D Laplace equation with potential
    perturbation. Top: numerical results. Bottom: cross-sections of
    $V(x)$ (left) and $\log_{10}|u(x)|$ (right) in the middle of the
    domain for the highest resolution test.}
  \label{fig:L31}
\end{figure}

\begin{figure}[h!]
  \begin{center}
    \begin{tabular}{|cc|cc|cc|}
      \hline
      $\max(|V|)$ & $N$ & $T_s$(sec) & $T_a$(sec) & $n_p$ & $T_p$(sec) \\
      \hline
      5.8e+02 & 1.2e+04 & 5.1e+00 & 3.8e-02 & 5 & 2.9e-01\\
      2.3e+03 & 1.0e+05 & 9.8e+01 & 2.6e-01 & 7 & 2.3e+00\\
      9.2e+03 & 8.6e+05 & 2.0e+03 & 2.8e+00 & 7 & 2.4e+01\\
      \hline
    \end{tabular}
    \includegraphics[height=1.8in]{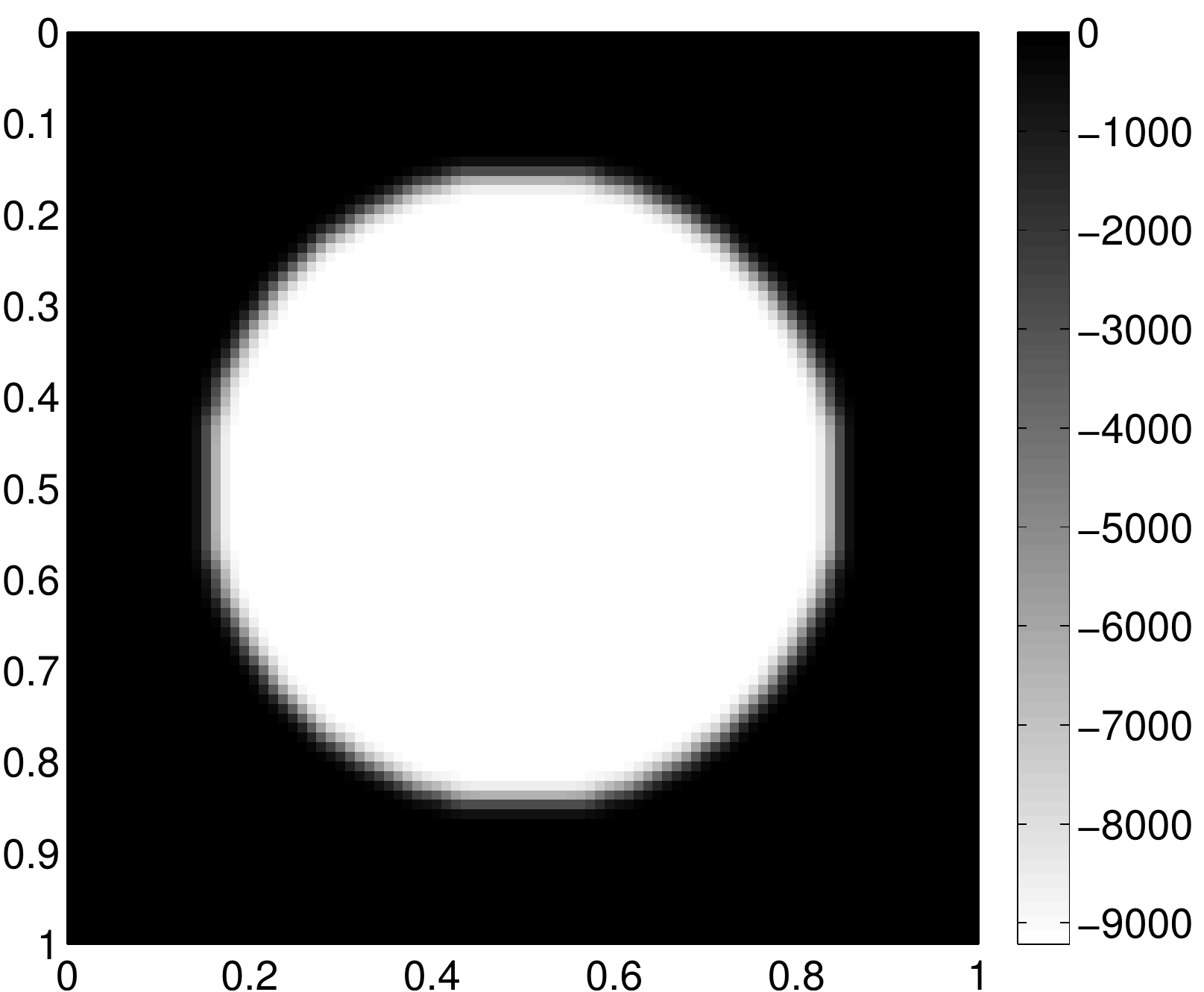} \hspace{0.25in} \includegraphics[height=1.8in]{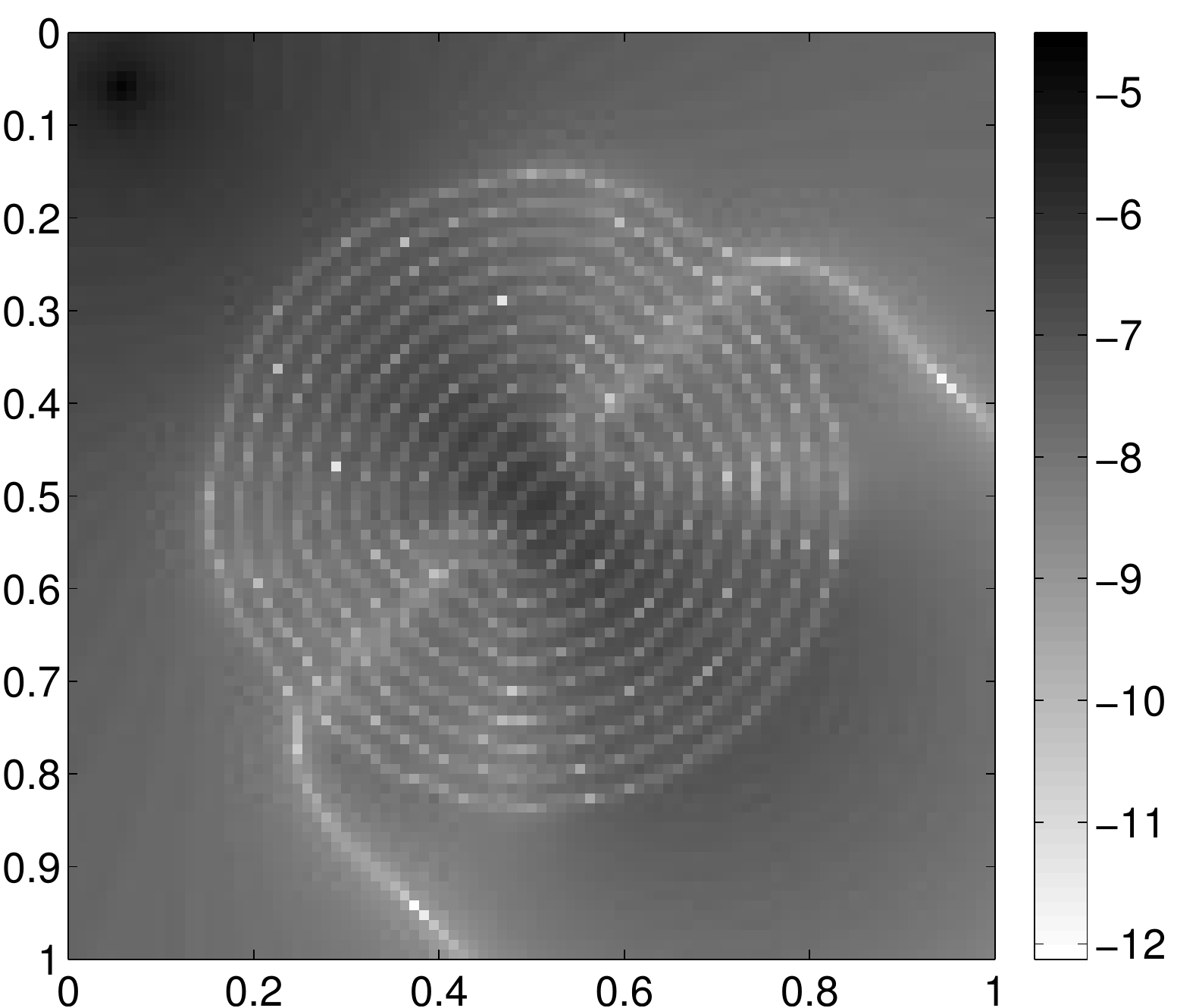}\\
  \end{center}
  \caption{Example 2 of 3D Laplace equation with potential
    perturbation. Top: numerical results. Bottom: cross-sections of
    $V(x)$ (left) and $\log_{10}|u(x)|$ (right) in the middle of the
    domain for the highest resolution test.}
  \label{fig:L32}
\end{figure}

The tests are performed for two examples of $V(x)$: a negative
Gaussian function and a smoothed characteristic function of an
$\ell_2$ ball. In each example, the right hand side $f$ is chosen to
be a delta source in $\Omega$ located near the top left corner of the
middle cross-section. The potential $V(x)$ is chosen to be negative
and proportional to $1/h^2$ in order to make sure that the problem
loses ellipticity. The numerical results are summarized in Figures
\ref{fig:L31} and \ref{fig:L32}.

The results demonstrate that the number of iterations grows at most
logarithmically with the problem size. This implies that the
numerically computed boundary condition operator $B$ is a good
approximation to the decaying condition of $u$ at infinity.

\section{Conclusion}

This paper introduces the sparsifying preconditioners for the
numerical solution of the Lippmann-Schwinger equation. The main idea
is to numerically transform this integral equation into a sparse and
localized linear system and then leverage existing efficient sparse
linear algebra algorithms. This approach combines the appealing
features of the integral equation formulation with the efficiency of
the PDE formulations. We discuss the algorithmic details for
rectangular and general domains, and extend this approach to the 3D
Laplace equation with a potential perturbation.

From the numerical results, we observe that the most of the
construction time of the preconditioners is spent on the construction
of the nested dissection factorization for the operator in
\eqref{eq:PC}. It is possible to replace this step with other solvers such as
the sweeping preconditioners \cite{engquist-2011a,engquist-2011b},
which are asymptotically more efficient.

Most of the discussion here is in the setting of inhomogeneous
acoustic scattering. An immediate task is to use this preconditioner
to study problems from electromagnetic and quantum scattering.

This approach is not restricted only to scattering in free space. A
future direction is to extend it to fully periodic systems or systems
that are periodic in certain directions. This should have direct
applications in the study of photonic crystal.

For most numerical methods, the stencil is defined analytically based
on the PDE. In the current approach, however, the stencils in the
preconditioning operators $A$ and $B$ are defined using the typical
behavior of the solution, either deterministically as in Section
\ref{sec:R} or randomly as in Section \ref{sec:G}. In general, such a
practice seems to be fairly unexplored in the field of numerical
analysis.

\bibliographystyle{abbrv} \bibliography{ref}

\end{document}